\documentclass[11pt]{amsart}
\usepackage{a4wide}
\usepackage{xcolor}
  \definecolor{darkred}{rgb}{0.4,0,0}
  \definecolor{darkgreen}{rgb}{0,0.4,0}
  \definecolor{darkblue}{rgb}{0,0,0.4}
\usepackage[colorlinks,citecolor=darkgreen,linkcolor=darkred,urlcolor=darkblue,filecolor=darkred]{hyperref}
\usepackage[noabbrev, nameinlink]{cleveref}
\usepackage{newpxtext}
\AtBeginDocument{%
  \DeclareFontShape{\encodingdefault}{\rmdefault}{m}{sl}{<-> pplro7t}{}%
  \DeclareFontShape{\encodingdefault}{\rmdefault}{b}{sl}{<-> pplbo7t}{}%
}
\usepackage{newpxmath}
\usepackage[style=alphabetic,giveninits=true,isbn=false,maxnames=99,maxalphanames=99]{biblatex}
\DeclareFieldFormat{title}{\mkbibquote{#1}}
\DeclareFieldFormat[misc]{date}{\iffieldequalstr{eprinttype}{arXiv}{\mkbibemph{Preprint} (#1)}{#1}}
\AtEveryBibitem{\ifentrytype{unpublished}{\clearfield{year}}{}}
\DeclareFieldFormat[unpublished]{note}{\mkbibemph{#1}}

\AtEveryBibitem{\clearfield{issn}}
\AtEveryCitekey{\clearfield{isbn}}
\AtEveryBibitem{\clearfield{number}}
\AtEveryBibitem{\ifentrytype{unpublished}{}{\ifentrytype{online}{}{\clearfield{url}}}}
\renewbibmacro{in:}{%
  \ifentrytype{article}{}{%
    \printtext{\bibstring{in}\intitlepunct}%
  }%
}
\usepackage{titlesec}
\titleformat{\section}{\normalfont\large\bfseries\scshape}{\mdseries\thesection.}{1em}{}
\titlespacing*{\section}{0pt}{.7\linespacing plus \linespacing}{.5\linespacing}
\titleformat{\subsection}[runin]{\bfseries}{\mdseries\thesubsection.}{.5em}{}[.]
\titlespacing*{\subsection}{0pt}{.5\linespacing plus .7\linespacing}{.5em}
\titleformat{\subsubsection}[runin]{\itshape}{\upshape\thesubsubsection.}{.5em}{}[.]
\titlespacing*{\subsubsection}{0pt}{.5\linespacing plus .7\linespacing}{.5em}
\makeatletter
\def\@settitle{\begin{center}%
  \baselineskip14\p@\relax
    \bfseries
    \large
\uppercasenonmath\@title
  \@title
  \end{center}%
}
\g@addto@macro\bfseries{\boldmath}
\makeatother

\usepackage{mathtools} 
\usepackage[shortlabels]{enumitem}
\usepackage{diagbox,booktabs}
\usepackage{caption}
\usepackage[labelfont=rm,labelformat=parens]{subcaption}

\usepackage{tikz}
\usetikzlibrary{decorations.pathreplacing,backgrounds}

\makeatletter
\tikzset{
    dot diameter/.store in=\dot@diameter,
    dot diameter=1pt,
    dot spacing/.store in=\dot@spacing,
    dot spacing=2pt,
    dots/.style={
        line width=\dot@diameter,
        line cap=round,
        dash pattern=on 0pt off \dot@spacing
    }
}
\makeatother

\newtheorem{theorem}{Theorem}[section]
\newtheorem{proposition}[theorem]{Proposition}
\newtheorem{corollary}[theorem]{Corollary}
\newtheorem{lemma}[theorem]{Lemma}
\theoremstyle{remark}
\newtheorem{remark}[theorem]{Remark}
\newtheorem*{notation}{Notation}
\providecommand{\defn}[1]{{\fontseries{b}\fontshape{sl}\selectfont #1}}
\makeatletter
\@definecounter{formula}
\let\c@formula\c@equation

\def\endformula{\eqno \hbox{\@eqnnum}$$\@ignoretrue}
\makeatother

\crefname{formula}{formula}{formulas}

\newcommand{\N}{\mathbf{N}}
\newcommand{\Z}{\mathbf{Z}}
\newcommand{\R}{\mathbf{R}}
\newcommand{\Hyp}{\mathbf{H}}
\newcommand{\PSL}{\operatorname{PSL}}
\newcommand{\pslz}{{\PSL(2,\Z)}}
\newcommand{\id}{{\mathrm{id}}}
\newcommand{\Id}{{\mathrm{Id}}}
\newcommand{\diag}{{\operatorname{diag}}}
\newcommand{\sphere}{\operatorname{S}}

\newcommand{\M}{\mathcal{M}}

\newcommand{\eval}{\operatorname{\pi}}
\newcommand{\cone}{\operatorname{\mathbf{C}}}
\newcommand{\suffix}{\operatorname{\mathbf{S}}}
\newcommand{\Gouezel}{\tilde{\suffix}}
\newcommand{\geodesics}{\operatorname{\mathbf{G}}}
\newcommand{\sufre}{\operatorname{\mathbf{E}}}
\newcommand{\typefunc}{\operatorname{\mathbf{T}}}

\begin{document}
\title[On the random walk on $\mathrm{PSL}(2,\mathbf{Z})$]
{Cone types and asymptotic invariants for the random walk on the modular group}
\author{Angel Pardo}
\thanks{This work was supported by ANID-Chile through the FONDECYT 1221934 grant.}
\address{
Departamento de Matem\'{a}tica y Ciencia de la Computaci\'{o}n,
Universidad de Santiago de Chile,
Las Sophoras 173, Estaci\'{o}n Central,
Santiago, Chile.
}
\email{angel.pardo@usach.cl}

\begin{abstract}
We compute the cone types of the Cayley graph of the modular group $\pslz$ associated with the standard system of generators {\small$\left\{\begin{psmallmatrix} 0 & -1 \\ 1 & 0 \end{psmallmatrix},\begin{psmallmatrix} 1 & 1 \\ 0 & 1 \end{psmallmatrix}\right\}$}.
We do this by showing that, in general, there is a set of suffixes of each element that completely determines the cone type of the element, and such suffixes are subwords of primitive relators.

Then, using J.~W.~Cannon's seminal ideas~(1984), we compute its growth function.
We estimate from above and below the spectral radius of the random walk using ideas from T.~Nagnibeda~(1999) and S.~Gouëzel~(2015).
Finally, using results of Y.~Guivarc'h (1980) and S.~Gouëzel, F.~Math\'{e}us and F.~Maucourant~(2015), we estimate other asymptotic invariants of the random walk, namely, the entropy and the drift.
\end{abstract}

\maketitle

\section{Introduction}

The modular group $\pslz$ is arguably one of the most fundamental groups in mathematics and, accordingly, one of the most studied.

In this note we study $\pslz$ through the lens of combinatorial and geometric group theory, focusing on the asymptotic invariants for the (simple symmetric) random walk on the Cayley graph of $\pslz$ associated with the standard generators $r=\begin{psmallmatrix} 0 & -1 \\ 1 & 0 \end{psmallmatrix}$ and $u=\begin{psmallmatrix} 1 & 1 \\ 0 & 1 \end{psmallmatrix}$. 

It is worth to mention that $\pslz$ is also the free product of the cyclic groups generated by $r$ and $ru=\begin{psmallmatrix} 0 & -1 \\ 1 & 1 \end{psmallmatrix}$, of order two and three, respectively.
Thus, we have the presentation $\pslz \cong \Z_2 * \Z_3 \cong \langle a,b | a^2,b^3\rangle$.
From the combinatorial and geometric group theoretic view point, this presentation is the handier.
In fact, the invariants for $\pslz$ associated with the generators $r$ and $ru$ can be computed explicitly, as we survey in the \Cref{sect:free-prod}.

We focus our attention in the generators $r$ and $u$ as they do not allow such straightforward computations but are still of interest. In fact, they are \emph{geometrically meaningful} as orbifold fundamental curves on the modular curve $\pslz\backslash\Hyp$ (as it is also the case for the generating set $\{ru,u\}$, which we address in \Cref{sect:other}).
Our case of study is also of interest from the combinatorial and geometric group theoretic point of view, since it is still a Cayley graph of a free product of two cyclic groups, but with respect to a generating set different from the ones that have been well understood.

We base our method on J.~W.~Cannon's classification of group elements by their cone types and our main contribution is the computation of the cone types of $\pslz$ relative to $\{r,u\}$.
This allows to compute the growth series using Cannon's original ideas~\cite{Cannon}, and to give numerical estimates for the spectral radius of the random walk using ideas from T.~Nagnibeda~\cite{Nagnibeda} and S.~Gouëzel~\cite{Gouezel}.
We also estimate the entropy and the drift of the random walk applying a result of S.~Gouëzel, F.~Math\'{e}us and F.~Maucourant~\cite{Gouezel-Matheus-Maucourant}.

One motivation for writing this note is that, even if it is likely that in the present case estimates (or even a complete description; cf. \Cref{sect:free-prod}) may be known by experts, we could not find any clue in the literature.
On the other hand, this work arises as a spin-off of \cite[Appendix~B]{Pardo:quantitative}, where we compute cone types for a different Fuchsian group in order to apply Nagnibeda's ideas and give lower bounds for the bottom of the spectrum of the combinatorial Laplace operator, which is equivalent to estimates from above for the spectral radius of the corresponding random walk, as shows \cref{form:bottom-radius} below.
In particular, there is some text overlap in respect to some background material.
We also stress that the estimates in the present work do not intend to be optimal in any sense.

\subsection{Cone types}
For a set $A$, we denote by $A^*$ the free monoid on $A$, an element in $A^*$ is called a \defn{word} in $A$, and the elements in $A$ are the \defn{letters}.

For a system of generators $S$ of a group $G$, in a slight abuse of notation we still write $S^*$ for $(S \cup S^{-1})^*$.
There is a natural evaluation morphism $\eval \colon S^* \to G$, that associates to a word $w \in S^*$ the product of its letters in $G$.
We denote the length of an element of $S^*$ by $|\cdot|$, that is, the number of its letters.
Similarly, we denote by $|\cdot|_S$ the \defn{word length} in $G$ associated with $S$.
That is, for an element $g \in G$, $|g|_S = \min\{|w|\colon w \in \eval^{-1}(g)\}$.
We say that $w \in S^*$ is a \defn{geodesic} if $|w| = |\eval(w)|_S$.
We also denote inversion by an overbar, that is, if $g \in G$, then $\bar g = g^{-1} \in G$ and, if $w = w_1 \dots w_n \in S^*$, where $w_1,\dots,w_n \in S$, then $\bar w = w_n^{-1} \dots w_1^{-1} \in S^*$.

Given $g \in G$, the \defn{cone} of $g$ relative to $S$, denoted $\cone(g)$, is the set of $h \in G$ for which some geodesic from $\id$ to $gh$ passes through $g$, that is,
\[
\cone(g) = \{h \in G\colon \exists v,w \in S^*, \eval(v) = g, \eval(w) = h, |vw| = |gh|_S\}.
\]
For each element $g \in G$, its cone $\cone(g)$ naturally defines a rooted subgraph of the Cayley graph.
We say that two cones $\cone_1, \cone_2$ are equivalent (or, of the same type) if they are isomorphic as rooted graphs.
The \defn{cone type} of $g \in G$ (relative to $S$) is the equivalence class of $\cone(g)$, which, in a slight abuse of notation, we still denote by $\cone(g)$.

It follows from Cannon's work that if $G$ is a word-hyperbolic group ---as is the case of the modular group $\pslz$---, then there are only finitely many cone types~\cite[Corollary~2]{Cannon}. 
Our proof follows the same lines but does not depend on this result, however.

\begin{theorem}[Cone types]
\label{thm:intro:cone-types}
Let $G=\pslz$ and $S = \{r,u\}$.
Then, there are exactly six cone types:
$\cone(\id)$, $\cone(r)$, $\cone(u)$, $\cone(ru)$, $\cone(ur)$ and $\cone(rur)$.
\end{theorem}

Furthermore, we give a simple combinatorial description for each cone type in \Cref{sect:cone-types}.
On the other hand, in the proof of \Cref{thm:intro:cone-types} arises the following general phenomenon of independent interest.

\begin{theorem}
\label{thm:intro:primitively-finite}
Let $G$ be a group generated by a finite system $S \subset G$ and suppose that the set of primitive relators (over $S$) is finite.
Then, $G$ has finitely many cone types (relative to $S$).
\end{theorem}

Here, by \defn{relator} we mean a word in $S$ that evaluates to the identity and by \defn{primitive}, that it does not contain proper subwords that are also relators.

In fact, we show in \Cref{prop:sufre-type}, that for any group $G$ and generating system $S \subset G$, the set of all maximal suffixes of geodesics for $g \in G$ that are subwords of a primitive or trivial relator completely determines the cone type of $g$.
In particular, knowing the set of primitive relators allows to compute, algorithmically at least, every cone type.

\subsection{Growth}
For $n \in \N$, let $\sphere_n \subset G$ be the sphere of radius $n$, that is, $\sphere_n = \{g \in G\colon |g|_S = n\}$.
The (spherical) \defn{growth series} of $G$ relative to $S$ is then the formal series
\[
V(x) = \sum_{g \in G} x^{|g|} = \sum_{n = 0}^{\infty} |\sphere_n|\, x^n.
\]
It follows from Cannon's work that if $G$ has finitely many cone types relative to $S$, then the growth series corresponds to a rational analytic function~\cite[Theorem~7]{Cannon} and the reciprocal of its radius of convergence gives the \defn{rate of exponential growth} of the group $G$ relative to $S$, denoted $v(G,S)$ and also called \defn{growth}, for short.

\begin{theorem}[Growth]
\label{thm:growth}
Let $G=\pslz$ and $S = \{r,u\}$.
Then, the growth series of $G$ relative to $S$ corresponds to the rational analytic function
\[
\frac{(1+x)(1+x+x^2)}{1-x-x^2}.
\]
In particular, the rate of exponential growth of $G$ relative to $S$ is the golden ratio, that is,
\[
v(G,S) = \varphi = \frac{1 + \sqrt{5}}{2} \approx
1.61 803.
\]
\end{theorem}

\begin{remark}
In particular, the growth sequence $(|\sphere_n|)_{n \in \N}$ corresponds to \cite[\href{https://oeis.org/A054886}{A054886}]{oeis}.
\end{remark}

\subsection{Random walk}
The (right, simple symmetric) \defn{random walk} on $G$ relative to $S$ is the Markov chain on $G$ whose transition probabilities are defined by
\[
p(g,h) = \begin{cases} \frac{1}{|S\cup \bar S|} & \text{ if } \bar gh \in S \cup \bar S, \\ \hfil 0 & \text{ otherwise}.\end{cases}
\]
A realization of the random walk starting from the identity is given by $X_0 = \id \in G$ and $X_n = s_1\dots s_n$, where $(s_i)_i$ is an independent sequence of $(S \cup \bar S)$-valued uniformly distributed random variables.
In other words, it is the (simple symmetric) random walk on the Cayley graph of $G$ relative to $S$.

\subsubsection{Spectral radius}
We denote $\M_S\colon \ell^2(G)\to\ell^2(G)$, the corresponding \defn{Markov operator}, that is,
\[
(\M_S h)_g \coloneqq \frac{1}{|S \cup \bar S|} \sum_{s\in S \cup \bar S} h_{gs},
\qquad
h\in\ell^2(G).
\]
We denote $\rho(G,S)$ the \defn{spectral radius} of the random walk on $G$ relative to $S$, that is, the spectral radius of $\M_S$,
\[
\rho(G,S) = \sup \left\{\left|\frac{\langle \M_S h, h \rangle}{\langle h, h \rangle}\right|,\; h\in \ell^2(G)\right\}.
\]
We estimate $\rho(G,S)$ from above following ideas of Nagnibeda~\cite{Nagnibeda}.
We also follow ideas of Gouëzel~\cite{Gouezel} to estimate $\rho(G,S)$ from below.
More precisely, we prove the following.

\begin{theorem}[Spectral radius]
\label{thm:spectral-radius}
Let $G=\pslz$ and $S = \{r,u\}$.
Then, the spectral radius of the random walk on $G$ relative to $S$ satisfies
\[
0.976 336
< \rho(G,S) <
0.976 642
.
\]
\end{theorem}

\begin{remark}
\label{rem:spectral-radius}
The spectral radius associated with a symmetric finite system of $d>1$ generators, is bounded from below by the spectral radius of (the random walk on) a regular tree of degree~$d$, that is, by $2\sqrt{d-1}/d$ (see, e.g., \cite[Chapter~II, Section~7.2]{ColinDeVerdiere}).
In our case, this yields
$
\rho(G,S) \geq 2\sqrt{2}/3 \approx
0.942 809
$.
In particular, our estimates from below improves the trivial bound obtained by comparison to the regular tree.

It is also worth to mention that the spectral radius of $G$ relative to $S' = \{r,ru\}$ can be explicitly computed.
In fact, 
$
\rho(G,S') = \frac{1}{6} + \frac{1}{6} \sqrt{13 + 8 \sqrt{2}} \approx
0.988 482
$
(see \Cref{free-prod:spectrum}).
\end{remark}

\subsubsection{Bottom of the spectrum of the Laplace operator}
A related quantity of interest is the \defn{bottom of the Laplace spectrum}, which we denote $\mu_0(G,S)$.
More precisely, let $\Delta_S\colon \ell^2(G)\to\ell^2(G)$ be the \defn{Laplace operator} on $G$ relative to $S$, that is,
\[
(\Delta_S h)_g \coloneqq \sum_{s\in S\cup\bar S} (h_g - h_{gs}),
\qquad
h\in\ell^2(G).
\]
Then $\mu_0(G,S)$ is the bottom of the spectrum of $\Delta_S$, that is,
\[\mu_0(G,S) = \inf \left\{\frac{\langle \Delta_S h, h \rangle}{\langle h, h \rangle},\; h\in \ell^2(G)\right\}.\]

By definition, $\Delta_S = |S \cup \bar S|(\Id_{\ell^2(G)} - \M_S)$.
Moreover, when the Cayley graph is bipartite ---as is the case for $G=\pslz$ and $S = \{r,u\}$---, the spectrum of the Markov operator is symmetric.
In such case, it follows that the bottom of the spectrum $\mu_0(G,S)$ is related to the spectral radius $\rho(G,S)$ by the formula 
\begin{formula}
\label{form:bottom-radius}
\mu_0(G,S) = |S \cup S^{-1}| (1-\rho(G,S)).
\end{formula}

As a direct consequence of \cref{form:bottom-radius} and \Cref{thm:spectral-radius} we get the following.
\begin{corollary}[Bottom of the Laplace spectrum]
\label{coro:combinatorial-spectrum}
Let $G=\pslz$ and $S = \{r,u\}$.
Then, the bottom of the spectrum of the Laplace operator on $G$ relative to $S$ satisfies
\[
0.0 700 754
< \mu_0(G,S) <
0.0 709 903
.
\pushQED{\qed}\qedhere
\]
\end{corollary}

\begin{remark}
Compare with the trivial upper bound given by the regular tree of degree three, that is,
$
\mu_0(G,S) \leq 3 - 2\sqrt{2} \approx
0. 171 572
$.
Also, compare with the bottom of the Laplace spectrum relative to $S' = \{r, ru\}$: 
$
\mu_0(G,S') = \frac{5}{2} - \frac{1}{2} \sqrt{13 + 8 \sqrt{2}} \approx
0.0 345 534
$
(see \Cref{free-prod:Laplace-spectrum}).
\end{remark}

\subsubsection{Entropy and drift}
Several numerical quantities have been introduced to describe the asymptotic behavior of random walks on groups (relative to a generating system). The spectral radius $\rho(G,S)$ is one of them.
Other important asymptotic invariants are the (asymptotic) \defn{entropy} $h(G,S)$ and the \defn{drift} $\ell(G,S)$ of the random walk. These are defined by
\[
h(G,S) = \lim_{n\to\infty} -\frac{1}{n} \sum_{g \in G} \upsilon_S^{*n}(g) \log(\upsilon_S^{*n}(g))
\qquad\text{and}\qquad
\ell(G,S) = \lim_{n\to\infty} \frac{1}{n} \sum_{g \in G} |g| \upsilon_S^{*n}(g),
\]
where $\upsilon_S$ is the uniform distribution on $S \cup \bar S$ and $\upsilon_S^{*n}$, its $n$-fold convolution.
There is a \emph{fundamental inequality} relating the entropy, the drift and the growth due to Guivarc'h~\cite{Guivarch}.
Moreover, Gouëzel, Math\'{e}us and Maucourant~\cite{Gouezel-Matheus-Maucourant} showed that estimates on the asymptotic invariants can be obtained from the others in several ways.
This, together with \Cref{thm:growth,thm:spectral-radius}, allows us to prove the following.

\begin{theorem}[Entropy and drift]
\label{thm:entropy-drift}
Let $G=\pslz$ and $S = \{r,u\}$.
Then, the entropy and the drift of the random walk on $G$ relative to $S$ satisfy
\[
0.0 938 046
< h(G,S) <
0.347 676
\qquad\text{and}\qquad
0.0 579 744
< \ell(G,S) <
0.214 876
.
\]
\end{theorem}

\subsection{Comparison between the three generating systems}
In \Cref{sect:free-prod,sect:other} we study other two \emph{geometrically meaningful} systems of generators of the modular group $\pslz$.
Namely, $S' = \{r,t\}$ and $S'' = \{t,u\}$, where $r = \begin{psmallmatrix} 0 & -1 \\ 1 & 0 \end{psmallmatrix}$, $t = \begin{psmallmatrix} 0 & -1 \\ 1 & 1 \end{psmallmatrix}$ and $u = \begin{psmallmatrix} 1 & 1 \\ 0 & 1 \end{psmallmatrix}$.
Note that $t = ru$.
For a schematic comparison, we summarize in \Cref{tab:comparative} the results on the asymptotic invariants associated with these three generating systems.

\begin{table}
\small
\begin{tabular}{ c c c c }
\toprule
\diagbox[font=\Small]{Invariant}{Generators} 
& $S = \{r,u\}$ & $S' = \{r,t\}$ & $S'' = \{t,u\}$
\\\midrule
Growth & $v = \frac{1+ \sqrt{5}}{2} \approx 1.618$ & $v = \sqrt{2} \approx 1.414$ & $v = 1+ \sqrt{2} \approx 2.414$
\\[1em]
Spectral radius & $0.9764 < \rho < 0.9767$ & $\rho = \frac{1+ \sqrt{13 + 8\sqrt{2}}}{6} \approx 0.9885$ & $0.8660 < \rho < 0.9268$
\\[1em]
Entropy & $0.09380 < h < 0.3477$ & $h = \frac{2\sqrt{2}}{15} \approx 0.1886$ & $0.2967 < h < 0.9069$
\\[1em]
Drift & $0.05797 < \ell < 0.2149$ & $\ell = \frac{2}{15} \approx 0.1333$ & $0.1229 < \ell < 0.3757$
\\\bottomrule
\end{tabular}
\caption{Asymptotic invariants associated to three systems of generators.}
\label{tab:comparative}
\end{table}

\subsection{Structure of the paper}
In  \Cref{sect:background}, we recall the basics of combinatorial group theory we need and prove some useful general combinatorial results. In particular, in \Cref{sect:general-cone-types} we prove \Cref{thm:intro:primitively-finite}.

In \Cref{sect:cone-types}, we compute the cone types of $G = \pslz$ relative to $S = \{r,u\}$, proving \Cref{thm:intro:cone-types}, and use them to produce a recurrence for the growth series as in Cannon's seminal work.
We prove \Cref{thm:cone-types} in \Cref{sect:cone-types} and \Cref{thm:growth}, in \Cref{sect:growth}.
We also include a graphical representation of the Cayley graph associated with $S$ in \Cref{fig:Cayley-graph}.

In \Cref{sect:spectral-radius} we estimate the spectral radius $\rho(G,S)$.
In \Cref{sect:Nagnibeda}, we state Nagnibeda's ideas for the upper bounds and, in \Cref{sect:Gouezel}, those of Gouëzel, for the lower bounds.
To conclude the proof of \Cref{thm:spectral-radius}, we include some details of the computations for the lower bound in \Cref{sect:suffix-type,sect:my-type}.

In \Cref{sect:entropy-drift} we state Guivarc'h~\cite{Guivarch} and Gouëzel--Math\'{e}us--Maucourant~\cite{Gouezel-Matheus-Maucourant} results relating the entropy and the drift to the growth and the spectral radius. Using our estimates for the spectral radius $\rho(G,S)$ and the exact value of the growth $v(G,S)$, this allows us to estimate the entropy $h(G,S)$ and the drift $\ell(G,S)$ and prove \Cref{thm:entropy-drift}.

In \Cref{sect:free-prod} we survey the analogous results in the case of the system of generators $S' = \{r,t=ru\}$, which yields a presentation of $\pslz$ as a free product of cyclic groups. In such case much more can be proven. In particular, we present the full Markov spectrum for the simple symmetric random walk in $\pslz$ relative to $S'$ and the exact values of the corresponding growth, spectral radius, entropy and drift.
We also include a graphical representation of the Cayley graph associated with $S'$ in \Cref{fig:Cayley-graph-1}.

In \Cref{sect:other} we repeat the study done for $S = \{r,u\}$, for the system of generators $S'' = \{t=ru,u\}$.
That is, we compute the growth series, estimate the spectral radius, the entropy and the drift of $\pslz$ relative to $S''$.
We also include a graphical representation of the Cayley graph associated with $S''$ in \Cref{fig:Cayley-graph-2}.

\subsection*{Acknowledgements}
The author is grateful Sebastien Gouëzel for the reference to the work of T. Nagnibeda.
The author is also greatly indebted to Tatiana Nagnibeda for useful comments and suggestions that improved the presentation of the paper, and for her patience providing him insight on the problems discussed in this work.

\section{Combinatorial group theory}
\label{sect:background}
In this section, we recall the basics of combinatorial group theory we need and prove some useful general combinatorial results useful for our purposes.

The following discussion is completely general.
For a complete introduction to this topic we refer the reader to Magnus--Karrass--Solitar's book~\cite{Magnus-Karrass-Solitar}.

Let $G$ be any group, and let $S$ be a subset of $G$. A \defn{word} in $S$ is any expression of the form
\[ w=s_{1}^{\sigma _{1}}s_{2}^{\sigma _{2}}\dots s_{n}^{\sigma _{n}}\]
where $s_1,\dots,s_n\in S$ and $\sigma_i \in \{+1,-1\}$, $i=1,\dots,n$. The number $|w|=n$ is the \defn{length} of the word.
The \defn{empy word}, denoted $\varepsilonup$, is the only one word of length zero.
Given two words $v$ and $w$, we say that $v$ is a \defn{subword} of $w$ if $w=xvy$, for some words $x$ and $y$. If $x$ is the empty word we say that $v$ is a \defn{prefix} of $w$. If $y$ is the empty word we say that $v$ is a \defn{suffix} of $w$.

Each word $w$ in $S$ represents an element of $G$.
Namely, the product of the expression, which we denote $\eval(w) \in G$.
For example, the identity element can be represented by the empty word, that is, $\id = \eval(\varepsilonup)$.
We say that two words are \defn{equivalent} if they represent the same element in $G$.

\begin{notation}
As usual, we use exponential notation for abbreviation; for example, the word $sss$ can also be denoted by $s^3$.
We also use an overbar to denote inverses, thus $\bar s$ stands for $s^{-1}$, 
and if $w = s_{1}^{\sigma _{1}}\dots s_{n}^{\sigma _{n}}$, then $\bar w = s_{n}^{-\sigma _{n}}\dots s_{1}^{-\sigma _{1}}$.
\end{notation}

In these terms, a subset $S$ of a group $G$ is a \defn{system of generators} if and only if every element of $G$ can be represented by a word in $S$.
Henceforth, let $S$ be a fixed system of generators of $G$ and a word is assumed to be a word in $S$.
We say that a non-empty word is a \defn{relator} if it represents the identity element of $G$, that is, if it equivalent to the empty word.
A generator next to its own inverse ($s\bar s$ or $\bar ss$) define a \defn{trivial relator}.

For a relator, we call a subword that is also a relator, a \defn{subrelator}.
We say that a word is \defn{reduced} if it has no trivial subrelators.
We say that a non-trivial relator is \defn{primitive} if if it does not contain proper subrelators.
In particular, a word is a geodesic if and only if it contains no primitive relators as subword.
Note that, if $S$ is a system of generators and $P$ is the set of all primitive relators on $S$, then $\langle S \mid P \rangle$ is a presentation of $G$.

For an element $g\in G$, we consider the \defn{word norm} $|g|_S$ to be the least length of a word which represents $g$ when considered as a product in $G$, and every such word is called a \defn{geodesic}, that is, if its length coincides with its word norm when considered as a product in $G$.
A geodesic does not contains relators as subwords and, in particular, it is is always reduced.
Also, a subword of a geodesic is also a geodesic.

The following decomposition result (see \Cref{figu:decomposition}) will be useful in \Cref{sect:cone-types} (in the form of \Cref{prop:primitive relator subword} below).

\begin{lemma}
\label{lemm:decomposition}
Let $v,w$ be two different equivalent geodesics. Then, there are geodesics $v_0,v_1$, $w_0,w_1$ and $x$ such that $v = v_0v_1x$ and $w = w_0w_1x$, and $v_1\bar w_1$ is a primitive relator (of even length).
\end{lemma}

\begin{figure}
\begin{tikzpicture}[scale=.2,blend mode=multiply]
\begin{scope}[line width=2,line cap=round]
\draw[blue] (4,0) ..controls (2,1.5) and (2,2).. (2,3) (2,3) ..controls (2,6) and (0,5).. (0,7) ..controls (0,8) and (0,8).. (-1,9) ..controls (-1.5,10) and (-.5,10).. (-.5,11) ..controls (-.5,11.5) and (-1,11.5).. (-1,12.5) -- (-1,14.5) ..controls (-1,16) and (1,16).. (1,18) ..controls (1,20) and (-1,19.5).. (-1,21) ..controls (-1,22) and (-1,23).. (1,24);

\draw[red] (4,0) ..controls (2,1.5) and (2,2).. (2,3) ..controls (2,5) and (-.5,4).. (-.5,5.5) ..controls (-.5,6) and (0,6).. (0,7) ..controls (0,8) and (0,8).. (-1,9) ..controls (-1.5,9.5) and (-1.75,10).. (-1.75,10.5) ..controls (-1.75,11.5) and (-1,11.5).. (-1,12.5) -- (-1,14.5) ..controls (-1,16) and (-2.5,16).. (-2.5,18) ..controls (-2.5,20) and (-1,19.5).. (-1,21) ..controls (-1,22) and (-1,23).. (1,24);
\end{scope}

\begin{scope}[line width=1]
\draw[blue!66!black,decoration={brace,mirror},decorate] (5,0+.1) -- node[midway,right=2] {$w_0$} (5,14.5-.1);
\draw[blue!66!black,decoration={brace,mirror},decorate] (5,14.5+.1) -- node[midway,right=2] {$w_1$} (5,21-.1);
\draw[blue!66!black,decoration={brace,mirror},decorate] (5,21+.1) -- node[midway,right=2] {$x$} (5,24-.1);
\draw[red!66!black,decoration={brace,mirror},decorate] (5,21+.1) -- node[midway,right=2] {$x$} (5,24-.1);

\draw[red!66!black,decoration={brace},decorate] (-4,0+.1) -- node[midway,left=2] {$v_0$} (-4,14.5-.1);
\draw[red!66!black,decoration={brace},decorate] (-4,14.5+.1) -- node[midway,left=2] {$v_1$} (-4,21-.1);

\draw[dots, very thin] (4.5,14.5) -- (-3.5,14.5);
\draw[dots, very thin] (4.5,21) -- (-3.5,21);
\end{scope}
\end{tikzpicture}
\caption{Decomposition of two equivalent geodesics.}
\label{figu:decomposition}
\end{figure}

\begin{proof}
Let $x$ be the largest common suffix of $v$ and $w$ (possibly $x$ is empty). Write $v=v'x$ and $w = w'x$.
Let $w_1$ and $v_1$ be the smallest non-empty suffixes of $w'$ and $v'$ respectively such that $v_1$ and $w_1$ are equivalent.
Such $v_1$ and $w_1$ exist since $v$ and $w$ are different words.
Moreover, they have the same length since they are equivalent geodesics, that is, they evaluate to the same element in $G$.
Write $v'=v_0v_1$ and $w'=w_0w_1$ (possibly $v_0$ and $w_0$ are empty).
In particular $v_0$ and $w_0$ are equivalent, since the same holds for $v',w'$ and $v_1,w_1$.
See \Cref{figu:decomposition}.

It remains to prove that $v_1\bar w_1$ is primitive.
Suppose $z$ is a subrelator of $v_1\bar w_1$.
Since $v_1$ and $w_1$ are geodesics, 
their subwords are also geodesics and, therefore, they do not contain any relator as subword.
It follows that $z = v_2\bar w_2$ for some non-empty suffixes $v_2$ and $w_2$, of $v_1$ and $w_1$, respectively. In particular, $v_2$ and $w_2$ are non-empty suffixes of $w'$ and $v'$, respectively, and $v_2,w_2$ are equivalent. But, by definition, $v_1$ and $w_1$ are the smallest such suffixes and therefore $v_2=v_1$ and $w_2=w_1$. Thus, $v_1\bar w_1$ has no proper subrelators and, therefore, $v_1\bar w_1$ is primitive.
\end{proof}

As a direct consequence of \Cref{lemm:decomposition}, we have the following.

\begin{proposition}
\label{prop:primitive relator subword}
Let $v=v'yx$ and $w=w'zx$ be two equivalent geodesics such that $y \bar z$ is also geodesic. Then, $y\bar z$ is a subword of some primitive relator (of even length).
\end{proposition}

\begin{proof}
Consider the decomposition given by \Cref{lemm:decomposition}. It is clear that $y$ is a subword of $v_1$ and $z$, of $w_1$. Then $y\bar z$ is a subword of the primitive relator $v_1\bar w_1$.
%
\end{proof}

\subsection{Cone types}
\label{sect:general-cone-types}
Recall that the \defn{cone} of an element $g \in G$, relative to $S$, is the set $\cone(g)$ of $h \in G$ for which some geodesic from $\id$ to $gh$ passes through $g$.
More precisely,  let $\geodesics(g)$ be the set of all geodesics for $g \in G$, that is, $\geodesics(g) = \{w \in S^*\colon \pi(w) = g, |w| = |g|_S\}$.
Then,
\[
\cone(g) = \{h \in G\colon \exists w \in \geodesics(g), v \in \geodesics(h) \colon wv \in \geodesics(gh)\}.
\]
Moreover, in a slight abuse of notation, we shall identify $\cone(g)$ with the rooted graph that it defines as subgraph of the Cayley graph of $G$.
Finally, the \defn{cone type} of an element $g \in G$ is the isomorphism class of rooted graphs of $\cone(g)$.

To any rooted graph $\Gamma$, one can associate its \defn{tree of geodesics}, that is, the (rooted) tree where vertices are geodesics from the root of $\Gamma$ and edges correspond to length-increasing edges in $\Gamma$.
It is clear that isomorphic rooted graphs yield isomorphic trees of geodesics.

In particular, in order to study cone types it is worth to determine their tree of geodesics.
A way to do this is through other \emph{type functions}.
A \defn{type function} for $G$ (relative to $S$) is a function $\tau \colon G \to T$, for some set $T$ of \defn{types}, such that the type of an element $\tau(g)$ determines the number of successors of each type.
In other words, two elements $g, g' \in G$ have the same type if and only if, for every $t\in T$,
\[
|\{s\in S_+(g) \colon \tau(gs)=t\}| = |\{s'\in S_+(g') \colon \tau(g's')=t\}|
,\]
where $S_\pm(g) = \{s \in S \cup \bar S \colon |gs| = |g| \pm 1\}$.

Thus, the tree of geodesics of cone types are, in a precise way, the minimal type functions.
In fact, it is straightforward from the definition that, for every $g \in G$, the type $\tau(g)$ completely determines the tree of geodesics of $\cone(g)$.

Consequently, once we have determined a type function, in order to completely determine the cone type of an element, one has to be able to detect geodesics representing the same vertex and possible length-preserving edges.
In the case where every relator has even length, there are no length-preserving edges and equivalent geodesics can be characterized using primitive relators through \Cref{prop:primitive relator subword}, for example.

\begin{remark}
\label{rem:fix-type}
In the presence of relators of odd length, it is convenient to be able to determine the type of neighbors of the same length as well.
Thus, an upgraded definition of type function in such case would include in addition the condition that, for every $g,g' \in G$ sharing the same type and for every type $t \in T$, $|\{s\in S_0(g) \colon \tau(gs)=t\}| = |\{s'\in S_0(g') \colon \tau(g's')=t\}|$, where $S_0(g) = \{s \in S \cup \bar S \colon |gs| = |g|\}$.
\end{remark}

The previous discussion motivates the following.

Let $\phi\colon S^* \to S^*$ be such that $\phi(w)$ is the longest suffix of $w \in S^*$ that is a subword of some primitive or trivial relator (we need to consider trivial relators in the presence of \emph{free letters}).

\begin{lemma}
\label{lemm:sufre}
Let $g,h \in G$.
Then, $h \in \cone(g)$ if and only if there is $v \in \geodesics(h)$ such that $\phi(w)v$ is a geodesic for every $w \in \geodesics(g)$.
\end{lemma}
\begin{proof}
Suppose that $h \in \cone(g)$, that is, there is $w_0 \in \geodesics(g)$ and $v \in \geodesics(h)$ such that $w_0v \in \geodesics(gh)$. In particular, $|gh|_S = |g|_S + |h|_S$.
Now, for every $w \in \geodesics(g)$, $|wv| = |w| + |v| = |g|_S + |h|_S = |gh|_S$ and therefore $wv \in \geodesics(gh)$.
Thus, $\phi(w)v$ is also a geodesic as subword of the geodesic $wv$.

Suppose now that $h \notin \cone(g)$, that is, $wv$ is not a geodesic for every $w \in \geodesics(g)$ and every $v \in \geodesics(h)$.
It follows that, for every $w \in \geodesics(g)$ and every $v \in \geodesics(h)$, there is a minimal non-empty suffix $w'$ of $w$ and a minimal non-empty prefix $v'$ of $v$ such that $w'v'$ is subword of a relator.
Now, choose $w'v'$ also minimal among all $w \in \geodesics(g)$ and $v \in \geodesics(h)$.
Then, similarly to the proof of \Cref{lemm:decomposition}, we can show that $w'v'$ is in fact subword of either a primitive or trivial relator.
It follows that $w'$ is a suffix of $\phi(w)$ and, therefore, $\phi(w)v$ is not a geodesic.
\end{proof}

Let $\sufre(g) = \phi(\geodesics(g))$, that is, $\sufre(g)$ is the set of all maximal suffixes of geodesics for $g$ that are subwords of a primitive or trivial relator.
Then, by \Cref{lemm:sufre}, we have the following.

\begin{proposition}
\label{prop:sufre-type}
The map $\sufre\colon G \to 2^{S^*}$ determines the cone types.
\qed
\end{proposition}

\begin{proof}[Proof of \Cref{thm:intro:primitively-finite}]
If there is a finite number of primitive relators, then $\sufre$ takes values in a finite set.
By \Cref{prop:sufre-type}, it follows that there is a finite number of cone types.
\end{proof}

Note that \Cref{lemm:sufre} not only shows that $\sufre$ determines the cone type of an element, but the cone itself.
Moreover, by construction, it is possible to determine the $\sufre$-type of the successors of a given $\sufre$-type.
In fact, we have the following.

\begin{lemma}
\label{lemm:successor-type}
Let $g \in G$. Then $s \in S_+(g)$ if and only if $s \in S_+(\eval(w))$ for every $w \in \sufre(g)$.
Moreover, for $s \in S_+(g)$, $\sufre(gs) = \sufre(\eval(\phi(ws)))$ for any $w \in \sufre(g)$.
\qed
\end{lemma}

\section{Cone types for the modular group}
\label{sect:cone-types}
The previous discussion is completely general.
We now specialize to the modular group $G=\pslz$ with generators $S = \{r,u\}$, where $r=\begin{psmallmatrix} 0 & -1 \\ 1 & 0 \end{psmallmatrix}$ and $u=\begin{psmallmatrix} 1 & 1 \\ 0 & 1 \end{psmallmatrix}$.

In order to compute the cone types, we start by computing $\sufre$, which is a type function by \Cref{prop:sufre-type}.
It is well known that $\langle r,u \mid r^2,(ru)^3\rangle$ is a presentation of $G$ (see, e.g., \cite{Alperin}).
Since we have the relator $r^2$, 
we can omit henceforth $\bar r$, as it coincides with $r$ as element in $G$.
The set of primitive relators is then 
given by
\[\{ r^2, (r u)^3, (r \bar u)^3, (ur)^3, (\bar ur)^3\}.\]

Denote $\suffix(g)$ the set of suffixes of geodesics for $g \in G$.
Then, by the description of the primitive relators, as a direct consequence of \Cref{prop:primitive relator subword}, we have the following.

\begin{corollary}
\label{coro:forbidden}
The following cases cannot happen:
\begin{align*}
\bullet \; & \; u,\bar u \in \suffix(g),
&
\bullet \; & \; ur,\bar ur \in \suffix(g),\\
\bullet \; & \; ar,\bar a^2 \in \suffix(g),
&
\bullet \; & \; ar,a \in \suffix(g),
&
\text{for } a & = u \text{ or } \bar u.
\end{align*}
\end{corollary}

\begin{proof}
Neither $u^2$, $ur\bar u$ nor $\bar uru$ are subwords of a primitive relator.
\end{proof}

It follows that $\sufre$ can only take values in the set
\[\{\emptyset,\{r\},\{a\},\{ra\},\{ar\},\{ra,a\},\{rar,\bar ar\bar a\}\}_{a=u,\bar u}\]
and any such value is possible.

\begin{proposition}
\label{prop:successors-type}
The $\sufre$-type of the successors of a given $\sufre$-type is determined as in \Cref{tab:diagram}.
\end{proposition}

\begin{table}
\begin{tabular}{@{\quad} >{$}l<{$} @{$\quad\to\quad$} >{$}l<{$} @{\quad}}
\toprule
\sufre(g) & \sufre(gs)\colon s\in S_+(g)\\
\midrule
\emptyset & \{r\}, \{u\}, \{\bar u\} \\
\{r\} & \{ru\}, \{r\bar u\} \\
\{a\} & \{ar\}, \{a\} \\
\{ra\} \text{ or } \{ra,a\} & \{a\}, \{rar,\bar ar\bar a\} \\
\{ar\} & \{ara,r\bar ar\}, \{r\bar a\} \\
\{rar,\bar ar\bar a\} & \{r\bar a,\bar a\} \\
\bottomrule
\end{tabular}
\caption{Each possible $\sufre$-type with the respective $\sufre$-type of its successors.
More precisely, for each $\sufre(g)$, $g\in G$, we show its respective multi-set of $\sufre(gs)$, $s\in S_+(g)$. Here, $a = u$ or $\bar u$.}
\label{tab:diagram}
\end{table}

\begin{proof}
Let $a\in \{u,\bar u\}$. Then, applying \Cref{lemm:successor-type}, we have the following:
\begin{itemize}
\item If $\sufre(g) = \emptyset$, then $g=\id$ and, evidently, $\sufre(gs) = \{s\}$, for $s \in \{r,u,\bar u\} = S_+(g)$.
\item If $\sufre(g) = \{r\}$, then $g=r$ and $\sufre(gb) = \{rb\}$, for $b\in\{u,\bar u\} = S_+(r)$.
\item If $\sufre(g) = \{a\}$, then $S_+(g)=\{r,a\}$, $\sufre(gr) = \{ar\}$ and $\sufre(ga) = \{a^2\}$.
\item If $\sufre(g) = \{ra\}$ or $\{ra,a\}$, then $S_+(g)=\{r,a\}$, $\sufre(gr) = \{ar,r\bar a\}$ and $\sufre(ga) = \{a^2\}$.
\item If $\sufre(g) = \{ar\}$, then $S_+(g)=\{u,\bar u\}$, $\sufre(ga) = \{ra,\bar ar\}$ and $\sufre(g\bar a) = \{r\bar a\}$.
\item If $\sufre(g) = \{rar,\bar ar\bar a\}$, then $S_+(g)=\{\bar a\}$ and $\sufre(g\bar a) = \{r\bar a,\bar a\}$.
\qedhere
\end{itemize}
\end{proof}

Recall that two cones $\cone_1, \cone_2$ are equivalent if they are isomorphic as rooted graphs.
And the \defn{cone type} of $g \in G$ (relative to $S$) is the equivalence class of $\cone(g)$, which, in a slight abuse of notation, we are still denoting by $\cone(g)$.

Now, even if the $\sufre$-types determine the cone types, the converse is not necessarily true.
In fact, it is clear from \Cref{tab:diagram} that there are different $\sufre$-types that \emph{behaves alike}.
Namely, the $\{ra\}$ and $\{ra,a\}$, as they share the multi-set of $\sufre$-types of their successors.
In our context, by \Cref{prop:sufre-type}, this implies that their cones are isomorphic.

Furthermore, it is clear that the map on $S^*$ that exchanges $u$ with $\bar u$ and fixes $r$ defines an automorphism of the Cayley graph of $\pslz$ relative to $S = \{r, u\}$, and consequently, several pairs with different $\sufre$-types define isomorphic cones though this map.
This motivates the definition of the following type function.
Let $\typefunc\colon G\to\{0,\dots,5\}$ be defined by
\[
\typefunc(g) =
\begin{cases}
0 & \text{if } \sufre(g) = \emptyset, \\
1 & \text{if } \sufre(g) = \{r\}, \\
2 & \text{if } \sufre(g) = \{u\} \text{ or } \{\bar u\}\\
3 & \text{if } \sufre(g) = \{ru\}, \{r\bar u\}, \{ru,u\} \text{ or } \{r\bar u,\bar u\}, \\
4 & \text{if } \sufre(g) = \{ur\} \text{ or } \{\bar ur\} \\
5 & \text{if } \sufre(g) = \{rur,\bar ur\bar u\} \text{ or } \{r\bar ur,uru\}.
\end{cases}
\]

\begin{theorem}
\label{thm:cone-types}
The function $\typefunc$ defines completely and uniquely the cone types.
Moreover,
\begin{itemize}
\item Type $0$ elements have one type $1$ and two type $2$ successors;
\item Type $1$ elements have two type $3$ successors;
\item Type $2$ elements have one type $2$ and one type $4$ successor;
\item Type $3$ elements have one type $2$ and one type $5$ successor;
\item Type $4$ elements have one type $3$ and one type $5$ successor; and
\item Type $5$ elements have one type $3$ successor.
\end{itemize}
\end{theorem}

\begin{proof}
From the previous discussion, it only remains to prove that different types define different cone types (cf. \Cref{fig:different-cone-types}).
We say that a vertex in a cone is at level $\ell \in \N$ if it is at distance $\ell$ from the root of the cone.
We say that a sub-cone is at level $\ell$ if its root is at level $\ell$.

It is clear that types $0$ and $5$ define different cone types from any other type, as they have respectively $3$ and $1$ vertices at level one, and all other types have two.
Similarly, types $3$ and $4$ define different cone types from types $1$ and $2$, as the former have a type $5$ cone at level one and the latter do not.

Now, type $1$ and $2$ define different cone types from each other, as the latter have a type $2$ cone at level one and the former have only type $3$ cones at level one. By the previous paragraph, we know that type $2$ and type $3$ cones are non-isomorphic.
Similarly, type $3$ and $4$ define different cone types from each other, as the former have type $2$ and $5$ cones at level one and the latter have type $3$ and $5$ cones at level one. Again, the claim holds since type $2$ and type $3$ cones are non-isomorphic.
\end{proof}

We include in \Cref{fig:different-cone-types} the first levels of the cones of each $\typefunc$-type from \Cref{thm:cone-types}, showing that they define different cone types.
We also include in \Cref{fig:cone-types} the \defn{graph of cone types}, where each vertex represents a cone type (labeled by the values given by $\typefunc$) and there is a directed edge from every cone type towards the cone type of each one of its successors.

\begin{figure}
\begin{subfigure}[b]{.618\textwidth}
\centering
\begin{tikzpicture}[node distance={15mm}, very thick, every node/.style = {draw, circle}]
\node (0) {$0$};
\node (1) [above right of=0] {$1$};
\node (2) [below right of=0] {$2$};
\node (3) [right of=1]{$3$};
\node (4) [right of=2]{$4$};
\node (5) [above right of=4] {$5$}; 

\begin{scope}[on background layer]
\begin{scope}[very thick]
\draw[->, red] (0) -- (1);
\draw[->,blue] (0) to [out=305,in=145,looseness=0] (2);
\draw[->,blue] (0) to [out=325,in=125,looseness=0] (2);

\draw[->,blue] (1) to [out=-10,in=190,looseness=0] (3);
\draw[->,blue] (1) to [out=10,in=170,looseness=0] (3);

\draw[->,blue] (2) to [out=270,in=180,looseness=4] (2);
\draw[->, red] (2) -- (4);

\draw[->,blue] (3) -- (2);
\draw[->, red] (3) to [out=325,in=125,looseness=0] (5);

\draw[->,blue] (4) -- (3);
\draw[->,blue] (4) -- (5);

\draw[->,blue] (5) to [out=145,in=305,looseness=0] (3);
\end{scope}
\end{scope}
\end{tikzpicture}
\vspace{-1em}
\caption{Graph of cone types.}
\label{fig:cone-types}
\end{subfigure}
\bigbreak
\begin{subfigure}[t]{.99\textwidth}
\centering
\begin{tikzpicture}[xscale=.4, very thick, every node/.style = {draw, circle, scale=.618}]
\node (0) at (0,0) {$0$};

\begin{scope}[on background layer]
\begin{scope}[very thick]
\draw[dots,red] (0) -- ++(0,.7);
\draw[dots,blue] (0) -- ++(-1,.7);
\draw[dots,blue] (0) -- ++(1,.7);
\end{scope}
\end{scope}
\end{tikzpicture}
\quad
\begin{tikzpicture}[xscale=.4, very thick, every node/.style = {draw, circle, scale=.618}]
\node (1) at (0,0) {$1$};
\node (13) at (2,1) {$3$};
\node (13') at (-2,1) {$3$};
\node (132) at (1,2) {$2$};
\node (135) at (3,2) {$5$};
\node (132') at (-3,2) {$2$};
\node (135') at (-1,2) {$5$};

\begin{scope}[on background layer]
\begin{scope}[very thick]
\draw[->,blue] (1) -- (13);
\draw[->,blue] (1) -- (13');
\draw[->,blue] (13) -- (132);
\draw[dots,blue] (132) -- ++(-.7,.7);
\draw[dots,red] (132) -- ++(.7,.7);
\draw[->,red] (13) -- (135);
\draw[dots,blue] (135) -- ++(0,.7);
\draw[->,blue] (13') -- (132');
\draw[dots,blue] (132') -- ++(-.7,.7);
\draw[dots,red] (132') -- ++(.7,.7);
\draw[->,red] (13') -- (135');
\draw[dots,blue]  (135')  -- ++(0,.7);
\end{scope}
\end{scope}
\end{tikzpicture}
\quad
\begin{tikzpicture}[xscale=.4, very thick, every node/.style = {draw, circle, scale=.618}]
\node (2) at (0,0) {$2$};
\node (22) at (-2,1) {$2$};
\node (24) at (2,1) {$4$};
\node (243) at (1,2) {$3$};
\node (245) at (3,2) {$5$};
\node (222) at (-3,2) {$2$};
\node (224) at (-1,2) {$4$};

\begin{scope}[on background layer]
\begin{scope}[very thick]
\draw[->,blue] (2) -- (22);
\draw[->,red] (2) -- (24);
\draw[->,blue] (22) -- (222);
\draw[dots,blue] (222) -- ++(-.7,.7);
\draw[dots,red] (222) -- ++(.7,.7);
\draw[->,red] (22) -- (224);
\draw[dots,blue] (224) -- ++(-.7,.7);
\draw[dots,blue] (224) -- ++(.7,.7);
\draw[->,blue] (24) -- (243);
\draw[dots,blue] (243) -- ++(-.7,.7);
\draw[dots,red] (243) -- ++(.7,.7);
\draw[->,blue] (24) -- (245);
\draw[dots,blue] (245) -- ++(0,.7);
\end{scope}
\end{scope}
\end{tikzpicture}
\quad
\begin{tikzpicture}[xscale=.4, very thick, every node/.style = {draw, circle, scale=.618}]
\node (3) at (0,0) {$3$};
\node (32) at (-2,1) {$2$};
\node (35) at (2,1) {$5$};
\node (322) at (-3,2) {$2$};
\node (324) at (-1,2) {$4$};

\begin{scope}[on background layer]
\begin{scope}[very thick]
\draw[->,blue] (3) -- (32);
\draw[->,red] (3) -- (35);
\draw[->,blue] (32) -- (322);
\draw[dots,blue] (322) -- ++(-.7,.7);
\draw[dots,red] (322) -- ++(.7,.7);
\draw[->,red] (32) -- (324);
\draw[dots,blue] (324) -- ++(-.7,.7);
\draw[dots,blue] (324) -- ++(.7,.7);
\draw[dots,blue] (35) -- ++(0,.7);
\end{scope}
\end{scope}
\end{tikzpicture}
\quad
\begin{tikzpicture}[xscale=.4, very thick, every node/.style = {draw, circle, scale=.618}]
\node (4) at (0,0) {$4$};
\node (43) at (-2,1) {$3$};
\node (45) at (2,1) {$5$};
\node (432) at (-3,2) {$2$};
\node (435) at (-1,2) {$5$};

\begin{scope}[on background layer]
\begin{scope}[very thick]
\draw[->,blue] (4) -- (43);
\draw[->,blue] (4) -- (45);
\draw[->,blue] (43) -- (432);
\draw[dots,blue] (432) -- ++(-.7,.7);
\draw[dots,red] (432) -- ++(.7,.7);
\draw[->,red] (43) -- (435);
\draw[dots,blue] (435) -- ++(0,.7);
\draw[dots,blue] (45) -- ++(0,.7);
\end{scope}
\end{scope}
\end{tikzpicture}
\quad
\begin{tikzpicture}[xscale=.4, very thick, every node/.style = {draw, circle, scale=.618}]
\node (5) at (0,0) {$5$};

\draw[dots,blue] (5) -- ++(0,.7);
\end{tikzpicture}
\caption{First levels of the cones defined by $\typefunc$.}
\label{fig:different-cone-types}
\end{subfigure}
\caption{Red arrows represent right multiplication by $r$ and blue, by $u$ or $\bar u$.}
\end{figure}

\subsection{Cayley graph and geometric description of cone types}
In the previous discussion, we have avoided natural geometric interpretations in terms of the Cayley graph.
The aim of this section is to give such a description for the ease of the reader.
However, this is not used elsewhere in this work.

Recall that the \defn{Cayley graph} $\Gamma(G,S)$ of a group $G$ relative to a system of generators $S \subset G$ is an edge-colored directed graph where the vertex set is $G$, the color set is $S \cup \bar S$.
If $g, h \in G$ and $s \in S \cup \bar S$, then there is a directed edge from $g$ to $h$ of color $s$ if and only if $h = gs$.
When there is $s \in S$ such that $s = \bar s$, we represent the two isochromatic opposite directed edges by an undirected one.

\begin{figure}
\includegraphics[width=.618\textwidth]{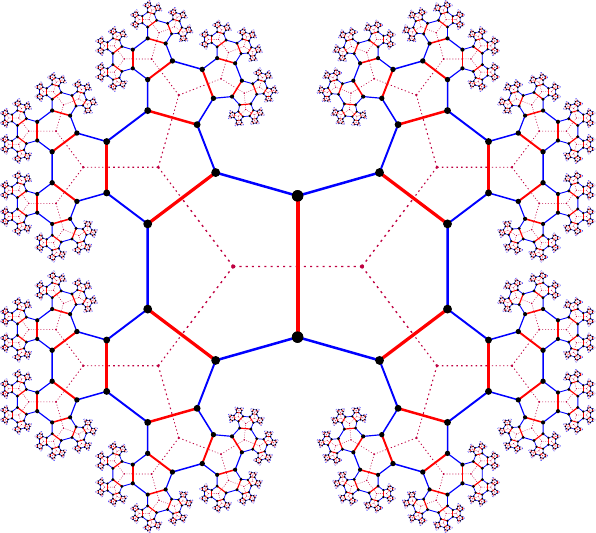}
\caption{The (undirected colored) Cayley graph $\Gamma(G,S)$ of $G = \pslz$ relative to $S = \{r,u\}$.
The $r$-edges in red and the $u$ and $\bar u$-edges in blue.
Also, in dashed purple lines, the weak dual regular tree of degree three.}
\label{fig:Cayley-graph}
\end{figure}

In the case of $G = \pslz$ and $S = \{r, u\}$, the Cayley graph $\Gamma(G,S)$ is weakly dual to a regular tree $\mathcal{T}$ of degree three, where the vertices in $\mathcal{T}$ correspond to hexagons ($6$-cycles) in $\Gamma(G,S)$ and the edges of $\mathcal{T}$ correspond to the $r$-edges of $\Gamma(G,S)$. See \Cref{fig:Cayley-graph} for a representation of the (undirected colored) Cayley graph $\Gamma(G,S)$.
Each red $r$-edge corresponds in fact to the $2$-cycle given by the relator $r^2$ and the hexagons, to the primitive relators of length six: $(ur)^3, (ru)^3, (\bar u)^3, (r \bar u)^3$.

Every element belong to two different hexagons and has exactly one neighbor (predecessor or successor) in each one of these different hexagons (corresponding to the $a$-edges, for $a = u$ and $\bar u$). We call these neigbors, the $u$-neighbors.
Similarly, every element has exactly one neighbor sharing the same two hexagons (corresponding to the $r$-edge), called $r$-neighbor.

On the other hand, each hexagon has a unique vertex of minimal length.
Thus, for a given $g \in G$, we can find a vertex $h \in G$ of minimal length, sharing a hexagon.
We say that the \defn{position} of $g$ is its distance to $h$, that is, $k = |\bar hg|_S = |g|_S - |h|_S \in \{0,1,2,3\}$.

The cone type of an element is then determined by its position and the combinatorics of the $r$-neighbor being a predecessor or a successor.
Note that, whether the $r$-neighbor is a predecessor or a successor is equivalent to consider the position of a point in the unique adjacent $2$-cycle determined by $r$.

The cone type ---as in \Cref{thm:cone-types}--- corresponds to the following combinatorics:
\begin{enumerate}[(1),start=0]
\item The position is $0$ and, in particular, the $r$-neighbor is a successor;
\item The position is $1$ and the $r$-neighbor is a predecessor;
\item The position is $1$ and the $r$-neighbor is a successor;
\item The position is $2$ and the $r$-neighbor is a predecessor;
\item The position is $2$ and the $r$-neighbor is a successor; and
\item The position is $3$ and, in particular, the $r$-neighbor is a predecessor.
\end{enumerate}

Reciprocally, this geometric description of the cone types can also be recovered from the combinatorial description given by the $\sufre$-types as follows.
The position of $g \in G$ is determined by the maximal length of an element in $\sufre(g)$ and whether or not there is $w \in \sufre(g)$ that ends in the letter $r$, determines whether the $r$-neighbor is a predecessor or a successor, respectively.

\subsection{Growth}
\label{sect:growth}
In this section we follow Cannon's ideas~\cite{Cannon} to compute the growth of $G = \pslz$ relative to $S = \{r,u\}$.

\begin{proof}[Proof of \Cref{thm:growth}]
Let $g$ be the generating function for the spherical growth sequence $a_n = |\sphere_n|$, $n \in \N$.
If $g_t$ is the generating function for the spherical growth sequence of elements of type $t \in \{0,\dots,5\}$, that is, of the sequence $a^{(t)}_n =  |\typefunc^{-1}(t) \cap \sphere_n|$, $n \in \N$.

Now, \Cref{thm:cone-types} (cf. \Cref{fig:cone-types}) implies that there is the following recurrence relation (cf. \cite[Proof of Theorem~7]{Cannon}):
\begin{gather*}
a^{(0)}_{0} = 1,
\qquad\text{and}\qquad
a^{(t)}_{0} = 0, \text{ for } t = 1,\dots,5
\\[1ex]
a^{(0)}_{n+1} = 0,
\qquad
a^{(1)}_{n+1} = a^{(0)}_{n},
\qquad
a^{(2)}_{n+1} = 2 a^{(0)}_{n} + a^{(2)}_{n} + a^{(3)}_{n}
\\[1ex]
a^{(3)}_{n+1} = 2 a^{(1)}_{n} + a^{(4)}_{n} + a^{(5)}_{n}
\qquad
a^{(4)}_{n+1} = a^{(2)}_{n},
\qquad
2 a^{(5)}_{n+1} = a^{(3)}_{n} + a^{(4)}_{n}
\end{gather*}
and, clearly, $a_{n} = a^{(0)}_{n} + a^{(1)}_{n} + a^{(2)}_{n} + a^{(3)}_{n} + a^{(4)}_{n} + a^{(5)}_{n}$.

It follows that the corresponding generating functions satisfy the equations
\begin{gather*}
g_0(z) = 1,
\qquad
g_1(z) = z g_0(z),
\qquad
g_2(z) = z (2 g_0(z) + g_2(z) +g_3(z)),
\\[1ex]
g_3(z) = z (2 g_1(z) + g_4(z) + g_5(z)),
\qquad
g_4(z) = z g_2(z),
\qquad
2 g_5(z) = z (g_3(z) + g_4(z))
\end{gather*}
and $g(z) = g_0(z) + g_1(z) + g_2(z) + g_3(z) + g_4(z) + g_5(z)$.

Finally, solving the recurrence, we get that
\begin{gather*}
g_0(z) = 1,
\qquad
g_1(z) = z,
\qquad
g_2(z) = \frac{2 z}{1 - z - z^2},
\\[1ex]
g_3(z) = \frac{2 z^2}{1 - z - z^2},
\qquad
g_4(z) = \frac{2 z^2}{1 - z - z^2},
\qquad
g_5(z) = \frac{2 z^3}{1 - z - z^2}
\end{gather*}
and \[g(z) = \frac{(1+z)(1+ z + z^2)}{1 - z - z^2}.\]

The growth rate $v(G,S)$ is the reciprocal of radius of convergence of $g$ around the origin, which is the root of smallest absolute value of  $z^2 + z - 1$. That is, it is the golden ratio
\[
v(G,S) = \varphi = \frac{1+\sqrt{5}}{2}.
\qedhere
\]
\end{proof}

\section{Spectral radius}
\label{sect:spectral-radius}
\subsection{Nagnibeda's ideas for the upper bound}
\label{sect:Nagnibeda}
In order to give upper bounds for the spectral radius $\rho(G,S)$, we follow Nagnibeda's ideas~\cite{Nagnibeda}, which are based in the following elementary result (cf. \cite[Chapter~II, Section~7.1]{ColinDeVerdiere}) and \cite[Section~1]{Nagnibeda}).

\begin{lemma}[Gabber--Galil]
\label{lemm:GG}
Let $G$ be a finitely generated group and $S$ a finite system of generators of $G$. Suppose there exists a function $L\colon G\times (S \cup \bar S)\to \R_+$ such that, for every $g\in G$ and $s\in S \cup \bar S$,
\[
L(g,s) = \frac{1}{L(gs,\bar{s})}
\qquad \text{ and } \qquad
\frac{1}{|S \cup \bar S|} \sum_{s\in S \cup \bar S} L(g,s) \leq \delta,
\]
for some $\delta>0$.
Then, $\rho(G,S) \leq \delta$.
\end{lemma}

\begin{proof}
Since
\[
\left(\sqrt{L(g,s)} h_{g} + \sqrt{L(gs,\bar{s})} h_{gs}\right)^2 \geq 0,
\]
we get that 
\[
2 h_{g} h_{gs}  \leq L(g,s) h_{g}^2 + L(gs,\bar{s}) h_{gs}^2.
\]
Summing over $g \in G$ and averaging over $s \in S \cup \bar S$, we get
\[
\langle \M_S h, h \rangle \leq \left(\frac{1}{|S \cup \bar S|} \sum_{s\in S \cup \bar S} L(g,s) \right) \langle h, h \rangle \leq \delta \langle h, h \rangle.
\]
That is, $\rho(G,S) \leq \delta$.
\end{proof}

For any type function $\tau\colon G \to T$ and positive valuation $c\colon T\to \R_+$, we can consider a function $L_c\colon G\times (S \cup \bar S)\to \R_+$ defined by
\[L_c(g,s) = \begin{cases} c_t, & \text{if } s \in S_+(g),\;t = \tau(gs), \\ 1/c_t, & \text{if } s\in S_-(g),\;t = \tau(g), \\ 1 & \text{otherwise}. \end{cases}\]
Then, every such $L_c$ satisfies $L_c(g,s) = 1/L_c(gs,\bar{s})$, since $s\in S_+(g)$ if and only if $\bar{s}\in S_-(gs)$.

Aditionally, for $t=\tau(g)\in T$, $g\in G$, we define
\[
f_t(c) \coloneqq \sum_{s\in S} L_c(g,s) = \sum_{s\in S_+(g)} c_{\tau(gs)} + \frac{1}{c_t}|S_-(g)| + |S| - |S_+(g)\cup S_-(g)|.
\]
This is well defined since $\tau$ is a type function and therefore the sum depends only on the type $t = \tau(g) \in T$.
Note that in the case where every relator has even length, we have that $S_+(g)\cup S_-(g) = S$, for every $g\in G$.

As a direct consequence of Gabber--Galil's \Cref{lemm:GG}, we get the following (cf. \cite[Section~2]{Nagnibeda}).

\begin{theorem}[Nagnibeda]
\label{thm:Nagnibeda}
Let $G$ be a finitely generated group and $S$ a finite system of generators of $G$.
Let $t\colon G\to \N$ be a compatible type function for $S$.
Then,
\[
\rho(G,S) \leq \frac{1}{|S \cup \bar S|} \sup_{k \in t(G)} f_k(c),
\]
for every $c\colon\N \to\R_+$, where $f_k$ is defined as above.
\qed
\end{theorem}

Then, every type function gives upper bounds for the spectral radius.

\begin{theorem}[Upper bound for the spectral radius]
\label{thm:upper-bound}
Let $G=\pslz$ and $S = \{r,u\}$.
Then, the spectral radius of the random walk on $G$ relative to $S$ satisfies
\[
\rho(G,S) <
0.976 642
.
\]
\end{theorem}
\begin{proof}
By \Cref{thm:cone-types}, the $f_t$'s of Nagnibeda's \Cref{thm:Nagnibeda} are given by:
\begin{align*}
f_0(c) & = c_1 + 2 c_2, &
f_1(c) & = 2 c_3 + 1/c_1, &
f_2(c) & = c_2 + c_4 + 1/c_2, \\
f_3(c) & = c_2 + c_5 + 1/c_3, &
f_4(c) & = c_3 + c_5 + 1/c_4, &
f_5(c) & = c_3 + 2/c_5.
\end{align*}
It follows that $\rho(G,S) \geq \max_{t} f_t(c)/|S \cup \bar S|$, for every $c=(c_1,\dots,c_5)\in \R_+^{5}$.
Thus, the problem can be reduced to find the optimal such bound.
This can be solved numerically: we get that $\bar c\in\R_+^5$ with $\bar c_1= 1,$
\[
\bar c_2 \approx 0.732 625 567, \quad
\bar c_3 \approx 0.792 704 707, \quad
\bar c_4 \approx 0.832 345 202, \quad
\bar c_5 \approx 0.935 795 167
\]
is a (local) minimun for $\max_{t} f_t(c)$, and
$\max_{t} f_t(\bar c) < 2.929 924 551$.

Finally, since $|S \cup \bar S|=3$, by Nagnibeda's \Cref{thm:Nagnibeda}, it follows that
\[
\rho(G,S) < 0.976 642.
\qedhere
\]
\end{proof}
\begin{remark}
The nature of Nagnibeda's estimates suggest that it is not possible to improve the upper bound using other type functions.
In fact, as shown by Nagnibeda~\cite[Section~3]{Nagnibeda:unbeatable}, these upper bounds correspond to the spectral radius of a random walk on the tree of geodesics of the group and, in particular, do not depend on the choice of the type function.
Additionally, it is also shown in \cite{Nagnibeda:unbeatable} that it is possible to compute this upper bound through a recurrence. However, we refrain from doing so here as the resulting estimates are the same (cf. \cite[Remark~3.3]{Nagnibeda:unbeatable} ).
\end{remark}

\subsection{Gouëzel's ideas for the lower bound}
\label{sect:Gouezel}
In order to give lower bounds for the spectral radius $\rho(G,S)$, we follow Gouëzel's ideas~\cite{Gouezel}.
The key tools are essentially the same type functions, but Gouëzel's techniques allows to neglect a finite number of elements of any type.

More precisely, following~\cite[Definition~1.2]{Gouezel}, we say that $\tau\colon G \to T$ is a \defn{type system} for $S$ if it is surjective, $T$ is finite and there is $M \in \N^{T \times T}$ such that, for all $i,j \in T$ and all but finitely many $g \in G$ with $\tau(g) = j$, we have
\[
|\{s\in S_+(g) \colon \tau(gs)=i\}| = M_{ij}.
\]

Then, Gouëzel's main result to estimate the spectral radius, \cite[Theorem~1.4]{Gouezel}, reads as follows.

\begin{theorem}[Gouëzel]
\label{thm:Gouezel}
Let $G$ be a group, finitely generated by $S \subset G$, without relators of odd length.
Let $\tau\colon G\to T$ be a type system for $S$ and suppose that the associated matrix $M $ is Perron--Frobenius.

Define the matrix $\tilde M$ by $\tilde M_{ij} = M_{ij} / p_i$, where $p_i$ is the number of predecessors of an element of type $i$, that is, $p_i = |S| - \sum_{j} M_{ij}$.
It follows that $\tilde M$ is also Perron--Frobenius and let $\eta > 0$ be its dominating eigenvalue, which is simple, and let $v \in \R_+^T$, an associated eigenvector with positive entries.

Let $D = \diag(v)$, $\hat M = D^{-1/2} M D^{1/2}$ and $\bar M = (\hat M + \hat M^{T})/2$.
Finally, let $\lambda$ be the maximal eigenvalue of the symmetric matrix $\bar M$.
Then,
\[
\rho(G,S) \geq \frac{2\lambda}{|S| \sqrt{\eta}}.
\]
\end{theorem}

\begin{theorem}[Lower bound for the spectral radius]
\label{thm:lower-bound}
Let $G=\pslz$ and $S = \{r,u\}$.
Then, the spectral radius of the random walk on $G$ relative to $S$ satisfies
\[
\rho(G,S) >
0.975 180
.
\]
\end{theorem}
\begin{proof}
Using the type function of \Cref{thm:cone-types}, by discarding all type $0$ and type $1$ elements, that is, $\id$ and $r$, respectively, we get a type system with $4 \times 4$ Perron-Frobenius matrices $M$ and $\tilde M$, where
\[
M = \begin{pmatrix} 1 & 1 & 0 & 0 \\ 0 & 0 & 1 & 1 \\ 1 & 0 & 0 & 0 \\ 0 & 1 & 1 & 0 \end{pmatrix}
\qquad\text{and}\qquad
\tilde M = \begin{pmatrix} 1 & 1 & 0 & 0 \\ 0 & 0 & 1 & 1 \\ 1 & 0 & 0 & 0 \\ 0 & 1/2 & 1/2 & 0 \end{pmatrix}.
\]
Following Gouëzel's \Cref{thm:Gouezel}, we compute the dominating eigenvalue $\eta$ and corresponding positive eigenvector $v \in \R_+^4$ of $\tilde M$, which are
\[
\eta = \frac{\sqrt{5}+1}{2}
\qquad\text{and}\qquad
v = \left( \frac{\sqrt{5}+1}{2},\, 1,\, 1,\, \frac{\sqrt{5}-1}{2} \right)^{T}.
\]
Then, we get $D = \diag(v)$, $\hat M = D^{-1/2} M D^{1/2}$ and
\[
\bar M = \frac{\hat M + \hat M^{T}}{2} = 
\frac{1}{2}
\begin{pmatrix}
2 & \sqrt{\frac{\sqrt{5}-1}{2}} & \sqrt{\frac{\sqrt{5}+1}{2}} & 0 \\[1ex]
\sqrt{\frac{\sqrt{5}-1}{2}} & 0 & 1 & \sqrt{\frac{\sqrt{5}+1}{2}}^3 \\[1ex]
\sqrt{\frac{\sqrt{5}+1}{2}} & 1 & 0 & \sqrt{\frac{\sqrt{5}+1}{2}} \\[1ex]
0 & \sqrt{\frac{\sqrt{5}+1}{2}}^3 & \sqrt{\frac{\sqrt{5}+1}{2}} & 0
\end{pmatrix}.
\]

Finally, we compute numerically the dominant eigenvalue of the symmetric matrix $\bar M$, which is
$\lambda > 1.860 673 779 029$.
Thus, by Gouëzel's \Cref{thm:Gouezel}, we get that
\[
\rho(G,S) \geq \frac{2\lambda}{3\sqrt{\eta}} > 0.975 180.
\qedhere
\]
\end{proof}

Note that this is not precisely the estimate from below given in \Cref{thm:spectral-radius}.
We present it, however to exhibit the computations in detail in a simpler case.
We conclude the proof of \Cref{thm:spectral-radius} in \Cref{sect:suffix-type,sect:my-type}, where we exhibit \emph{finer} type systems and give the corresponding numerical lower bounds using Gouëzel's \Cref{thm:Gouezel}.

\section{Entropy and drift}
\label{sect:entropy-drift}
In addition to the spectral radius, other important numerical quantities have been introduced to describe the asymptotic behavior of random walks on groups.
In this section we estimate some of the most relevant asymptotic invariants: the (asymptotic) \defn{entropy} $h(G,S)$ and the \defn{drift} $\ell(G,S)$ of the random walk. These are defined by
\[
h(G,S) = \lim_{n\to\infty} -\frac{1}{n} \sum_{g \in G} \upsilon_S^{*n}(g) \log(\upsilon_S^{*n}(g))
\qquad\text{and}\qquad
\ell(G,S) = \lim_{n\to\infty} \frac{1}{n} \sum_{g \in G} |g| \upsilon_S^{*n}(g),
\]
where $\upsilon_S$ is the uniform distribution on $S \cup \bar S$ and $\upsilon_S^{*n}$, its $n$-fold convolution.
Gouëzel, Math\'{e}us and Maucourant~\cite{Gouezel-Matheus-Maucourant} showed that estimates on these numerical quantities can be obtained from the others in several ways.
More precisely,
consider the strictly increasing function $F \colon [0,1) \to \R_+$ defined by 
\[
F(x) = x \log \left(\frac{1+x}{1-x}\right) = 2x \operatorname{arctanh}(x)
.\]
Then, for any (simple, symmetric) random walk on a group, they proved the following (cf. \cite[Theorem~1.1]{Gouezel-Matheus-Maucourant}).

\begin{theorem}[Gouëzel--Matheus--Maucourant]
\label{thm:GMM}
We have that $h \geq F\left(\sqrt{1-\rho^2}\right)$ and $h \geq F\left(\ell\right)$.
\end{theorem}

Previously, Guivarc'h~\cite{Guivarch} proved the so-called \emph{fundamental inequality} between the entropy $h$, the drift $\ell$ and the growth $v$.

\begin{theorem}[Guivarc'h]
\label{thm:Guivarc'h}
We have that $h \leq v \ell$.
\end{theorem}

\begin{proof}[Proof of \Cref{thm:entropy-drift}]
By Gouëzel--Matheus--Maucourant's \Cref{thm:GMM}, $h \geq F\left(\sqrt{1-\rho^2}\right)$ and, by \Cref{thm:spectral-radius}, $\rho < \bar \rho \coloneqq 0.976 641 504$.
It follows that 
\[
h > \bar h \coloneqq F\left(\sqrt{1-\bar\rho^2}\right) > 
0.0 938 046 6.
\]
We also have that $F(\ell) \leq h$. Since $F$ is increasing in $[0,1)$, we have that
\[
\ell \leq F^{-1}\left(h\right) < \bar\ell \coloneqq F^{-1}\left(\bar h\right) <
0.214 875 3.
\]

On the other hand, $v = \frac{1+\sqrt{5}}{2}$, by \Cref{thm:growth}.
Thus, by Guivarc'h's \Cref{thm:Guivarc'h}, we get
\[
h \leq v \bar\ell <
0.347 675 5
\qquad\text{and}\qquad
\ell \geq \bar h/v >
0.0 579 744 6.
\qedhere
\]
\end{proof}

\appendix
\section{The modular group as a free product of cyclic groups. A brief survey}
\label{sect:free-prod}

In this appendix, we include a brief survey on analogous results in the case of the generating system $S'=\{r,t\}$, where $r = \begin{psmallmatrix} 0 & -1 \\ 1 & 0 \end{psmallmatrix}$ and $t = \begin{psmallmatrix} 0 & -1 \\ 1 & 1 \end{psmallmatrix}$.
Note that, one has $t = ru$, for $u = \begin{psmallmatrix} 1 & 1 \\ 0 & 1 \end{psmallmatrix}$.

With this generating system, we have the presentation $\pslz = \langle r,t \mid r^2, t^3 \rangle$.
In particular, as abstract group, $\pslz$ is the free product of the cyclic groups of order two and three, that is, $\pslz \cong \Z_2 * \Z_3$ and the elements in $S'$ correspond to generators for the cyclic groups (see, e.g., \cite{Alperin}).

The Cayley graph of such a free product relative to the cyclic generators has a particularly easy structure and a much more detailed descriptions of the objects studied in this work is possible.
For example, one can easily derive the growth series.
In fact, it is enough to count the reduced words in $S'^*$ with no repeated consecutive letters.
A recurrence can then be easily derived by counting such words that end with an $r$ separately from those ending in $t$ or $\bar t$.

\begin{theorem}[Growth]
\label{thm:growth-1}
Let $G = \pslz$ and $S' = \{r,t\}$.
Then, the growth series of $G$ relative to $S'$ corresponds to the rational analytic function
\[
\frac{1+3z+2z^2}{1-2x^2}.
\]
In particular, the rate of exponential growth of $G$ relative to $S'$ is
\[
v(G,S') = \sqrt{2} \approx
1.414 21
.
\pushQED{\qed}\qedhere
\]
\end{theorem}

This is not different in essence from what we do for the generating system $S = \{r,u\}$ in \Cref{sect:cone-types}.
However, in this case, the simpler structure allows an easier understanding of the combinatorics.
In fact, geodesics words are exactly the reduced words with no consecutive occurrences of letters, and the ending letter of such a geodesic determines the cone type of the corresponding group element.

\begin{theorem}[Cone types]
\label{thm:cone-types-1}
Let $G=\pslz$ and $S' = \{r,t\}$.
Then, there are exactly three cone types:
$\cone(\id)$, $\cone(r)$ and $\cone(t)$.
\pushQED{\qed}\qedhere
\end{theorem}

With considerably more work, it is also possible to describe completely the spectrum of the simple symmetric random walk on $G = \pslz$ associated with $S' = \{r,t\}$ (see \cite[Theorem~4]{Gutkin}; originally proved in \cite{McLaughlin}).

\begin{theorem}[Markov spectrum]
\label{free-prod:spectrum}
Let $G=\pslz$ and $S' = \{r,t\}$. And consider the Markov operator $\M_{S'}$ associated with the simple symmetric random walk on $G$ relative to $S'$.
Then, the point spectrum of $\M_{S'}$ is $\left\{-\frac{2}{3},0\right\}$ and the absolutely continuous spectrum of $\M_{S'}$ is
\[
\left[\frac{1}{6} - \frac{1}{6} \sqrt{13 + 8 \sqrt{2}}, \frac{1}{6} - \frac{1}{6} \sqrt{13 - 8 \sqrt{2}}\right]
\cup
\left[\frac{1}{6} + \frac{1}{6} \sqrt{13 - 8 \sqrt{2}}, \frac{1}{6} + \frac{1}{6} \sqrt{13 + 8 \sqrt{2}}\right].
\]
In particular, the spectral radius of the random walk on $G$ relative to $S$ is
\[
\rho(G,S') = \frac{1}{6} + \frac{1}{6} \sqrt{13 + 8 \sqrt{2}} \approx
0.988 482
.
\pushQED{\qed}\qedhere
\]
\end{theorem}

\begin{corollary}[Laplace spectrum]
\label{free-prod:Laplace-spectrum}
Let $G=\pslz$ and $S' = \{r,t\}$. And consider the combinatorial Laplace operator $\Delta_{S'}$ on $G$ relative to $S'$.
Then, the point spectrum of $\Delta_{S'}$ is $\left\{3,5\right\}$ and the absolutely continuous spectrum of $\Delta_{S'}$ is
\[
\left[\frac{5}{2} - \frac{1}{2} \sqrt{13 + 8 \sqrt{2}}, \frac{5}{2} - \frac{1}{2} \sqrt{13 - 8 \sqrt{2}}\right]
\cup
\left[\frac{5}{2} + \frac{1}{2} \sqrt{13 - 8 \sqrt{2}}, \frac{5}{2} + \frac{1}{2} \sqrt{13 + 8 \sqrt{2}}\right].
\]
In particular, the bottom of the spectrum of the Laplace operator on $G$ relative to $S$ is
\[
\mu_0(G,S') = \frac{5}{2} - \frac{1}{2} \sqrt{13 + 8 \sqrt{2}} \approx
0.0345534
.
\pushQED{\qed}\qedhere
\]
\end{corollary}

Furthermore, it is also possible to compute the other asymptotic invariants of a random walks on a free product of cyclic groups (see \cite[Sections~4.2 and 5.1]{Mairesse-Matheus}).

\begin{theorem}[Asymptotic invariants]
\label{free-prod:invariants}
Let $G=\pslz$ and $S' = \{r,t\}$.
Then, the rate of exponential growth $v = v(G,S')$, the entropy $h = h(G,S')$ and the drift $\ell = \ell(G,S')$ satisfy $h = v \ell$.
Moreover,
\[
h = \frac{2\sqrt{2}}{15} \approx
0.188562
\qquad\text{and}\qquad
\ell = \frac{2}{15} \approx
0.133333
.
\pushQED{\qed}\qedhere
\]
\end{theorem}

\subsection*{Cayley graph}
\begin{figure}
\includegraphics[width=.618\textwidth]{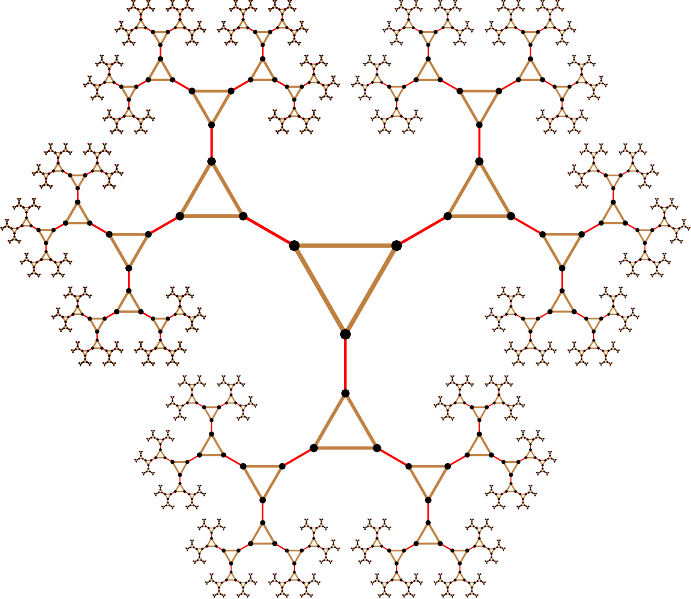}
\caption{The (undirected colored) Cayley graph $\Gamma(G,S')$ of $G = \pslz$ relative to $S' = \{r,t\}$.
The $r$-edges in red and the $t$ and $\bar t$-edges in brown.}
\label{fig:Cayley-graph-1}
\end{figure}

In the case of $G = \pslz$ and $S' = \{r, t\}$, the Cayley graph $\Gamma(G,S)$ has a minor regular tree $\mathcal{T}$ of degree three, where each vertex in $\mathcal{T}$ correspond to a triangle ($3$-cycles) in $\Gamma(G,S')$ and the edges of $\mathcal{T}$ correspond to the $r$-edges of $\Gamma(G,S')$. See \Cref{fig:Cayley-graph-1} for a representation of the (undirected colored) Cayley graph $\Gamma(G,S')$.
Each red $r$-edge corresponds in fact to the $2$-cycle given by the relator $r^2$ and the triangles, to the primitive relators of length three: $t^3$ and $\bar t^3$.

In particular, every element belong to exactly one triangle and every triangle contains a single vertex of minimal length, determining the cone type of any element other than the identity.
In fact, an element $g \in G\setminus\{\id\}$ is of minimal length in its triangle if and only if its only geodesic ends with the letter $r$.

\section{Another geometrically meaningful generating system}
\label{sect:other}
In this appendix, we include a summary of the analogous results in the case of the generating system $S''=\{t,u\}$, where $t = \begin{psmallmatrix} 0 & -1 \\ 1 & 1 \end{psmallmatrix}$ and $u = \begin{psmallmatrix} 1 & 1 \\ 0 & 1 \end{psmallmatrix}$.
Note that $t = ru$, where $r = \begin{psmallmatrix} 0 & -1 \\ 1 & 0 \end{psmallmatrix}$.

In this case, we have the presentation $\langle t,u \mid t^3, (t\bar u)^2 \rangle$ and the set of primitive relators is
\[
\{t^3, \bar t^3, (t\bar u)^2, (\bar tu)^2, (\bar ut)^2, (u\bar t)^2\}.
\]

\begin{theorem}[Cone types]
\label{thm:cone-types-2}
Let $G=\pslz$ and $S'' = \{t,u\}$.
Then, there are exactly three cone types for $G$ relative to $S''$:
$\cone(\id)$, $\cone(t)$ and $\cone(u)$.
\end{theorem}

\begin{proof}
Since all primitive relators are of length at most four, any suffix of a geodesic that is a subword of a primitive relator has length at most two.
It follows that the type function $\sufre$ from \Cref{prop:sufre-type} can only take values in the set
\[
\{\emptyset, \{t\}, \{\bar t\}, \{u\}, \{\bar u\}, \{t\bar u, u \bar t\}, \{\bar tu, \bar ut\}\}
\]
and any such value is possible.
In \Cref{tab:diagram-2}, we describe the $\sufre$-type of their successors.
Such description follows easily from the definition of $\sufre$ and the description of the primitive relators above.
We also include the types of the neighbors of the same length.
This is needed as there are relators of odd length (see \Cref{rem:fix-type}).

Now, it is clear also from the primitive relators that the involution that interchange each letter in $S''$ with its inverse, defines an automorphism of the corresponding Cayley graph and, in particular, isomorphisms between the respective cones. Thus, we have $\cone(t) \cong \cone(\bar t)$, $\cone(u) \cong \cone(\bar u)$, $\cone(\bar t\bar u) \cong \cone(tu)$ and $\cone(\bar ut) \cong \cone(u\bar t)$.

Moreover, it is still possible to reduce the number of cone types noticing that $\cone(ut) \cong \cone(t)$.
For simplicity, let $\typefunc''\colon G\to\{0,\dots,2\}$ be defined as follows:
\[
\typefunc''(g) =
\begin{cases}
0 & \text{if } \sufre(g) = \emptyset, \\
1 & \text{if } \sufre(g) = \{t\}, \{\bar t\}, \{t\bar u, u\bar t\} \text{ or } \{\bar tu, \bar ut\} \\
2 & \text{if } \sufre(g) = \{u\} \text{ or } \{\bar u\}.
\end{cases}
\]
Then, by the description of the $\sufre$-types in \Cref{tab:diagram-2}, it is clear that $\typefunc''$ is a type function.
In particular, $\typefunc''$ determines the cone types (see also \Cref{fig:cone-types-2}).

It only remains to prove that different $\typefunc''$-types define different cone types.
But this follows easily, similarly as in the proof of \Cref{thm:cone-types}, simply counting the number of elements in the level one in each rooted graph (see \Cref{fig:different-cone-types-2}).
\end{proof}

\begin{table}
\begin{tabular}{@{\quad} >{$}l<{$} @{$\quad\to\quad$} >{$}l<{$} @{\quad}}
\toprule
\sufre(g) & \sufre(gs)\colon s\in S''_+(g) \\
\midrule
\emptyset & \{t\}, \{\bar t\}, \{u\}, \{\bar u\} \\
\{t\} & \{t\bar u, u\bar t\}, \{u\} \\
\{\bar t\} & \{\bar tu, \bar ut\}, \{\bar u\} \\
\{u\} & \{t\bar u, u\bar t\}, \{t\}, \{u\} \\
\{\bar u\} & \{\bar tu,\bar ut\}, \{\bar t\}, \{\bar u\} \\
\{t\bar u, u\bar t\} & \{\bar t\}, \{\bar u\} \\
\{\bar tu, \bar ut\} & \{t\}, \{u\} \\
\bottomrule
\end{tabular}
\hfil
\begin{tabular}{@{\quad} >{$}l<{$} @{$\quad\to\quad$} >{$}l<{$} @{\quad}}
\toprule
\sufre(g) & \sufre(gs)\colon s\in S''_0(g) \\
\midrule
\emptyset & \\
\{t\} & \{t\}\\
\{\bar t\} & \{\bar t\} \\
\{u\} & \\
\{\bar u\} & \\
\{t\bar u, u\bar t\} & \{\bar t\} \\
\{\bar tu, \bar ut\} & \{t\} \\
\bottomrule
\end{tabular}
\caption{The possible $\sufre$-types relative to $S''$ with the $\sufre$-type of their respective successors (left) and neighbors of the same length (right).}
\label{tab:diagram-2}
\end{table}

\begin{figure}
\begin{subfigure}[t]{0.618\textwidth}
\centering
\begin{tikzpicture}[xscale=.5, very thick, every node/.style = {draw, circle, scale=.618}]
\node (0) at (0,0) {$0$};
\node (01) at (-3,1) {$1$};
\node (01') at (-1,1) {$1$};
\node (02) at (1,1) {$2$};
\node (02') at (3,1) {$2$};

\begin{scope}[on background layer]
\begin{scope}[very thick]
\draw[->,brown] (0) -- (01);
\draw[->,brown] (0) -- (01');

\draw[dash pattern=on 2 off 1,brown] (01) -- (01');
\draw[dots,blue] (01) -- ++(-.7,.7);
\draw[dots,blue] (01) -- ++(.7,.7);
\draw[dots,blue] (01') -- ++(-.7,.7);
\draw[dots,blue] (01') -- ++(.7,.7);

\draw[->,blue] (0) -- (02);
\draw[->,blue] (0) -- (02');

\draw[dots,brown,line cap=round] (02) -- ++(-.9,.7) ++(.2,.05) coordinate (*);
\draw[dash pattern=on 2 off 1,brown,semithick] (*) -- ++(.6,0);
\draw[dots,brown] (02) -- ++(0,.7);
\draw[dots,blue] (02) -- ++(.9,.7);
\draw[dots,brown] (02') -- ++(-.9,.7) ++(.2,.05) coordinate (*);
\draw[dash pattern=on 2 off 1,brown,semithick] (*) -- ++(.6,0);
\draw[dots,brown] (02') -- ++(0,.7);
\draw[dots,blue] (02') -- ++(.9,.7);
\end{scope}
\end{scope}
\end{tikzpicture}
\quad
\begin{tikzpicture}[xscale=.5, very thick, every node/.style = {draw, circle, scale=.618}]
\node (1) at (0,0) {$1$};
\node (11) at (-1,1) {$1$};
\node (12) at (1,1) {$2$};

\begin{scope}[on background layer]
\begin{scope}[very thick]

\draw[->,blue] (1) -- (11);
\draw[->,blue] (1) -- (12);

\draw[dots,blue] (11) -- ++(-.7,.7);
\draw[dots,blue] (11) -- ++(.7,.7);
\draw[dots,brown] (12) -- ++(-.9,.7) ++(.2,.05) coordinate (*);
\draw[dash pattern=on 2 off 1,brown,semithick] (*) -- ++(.6,0);
\draw[dots,brown] (12) -- ++(0,.7);
\draw[dots,blue] (12) -- ++(.9,.7);
\end{scope}
\end{scope}
\end{tikzpicture}
\quad
\begin{tikzpicture}[xscale=.5, very thick, every node/.style = {draw, circle, scale=.618}]
\node (2) at (0,0) {$2$};
\node (21) at (-2,1) {$1$};
\node (21') at (0,1) {$1$};
\node (22) at (2,1) {$2$};

\begin{scope}[on background layer]
\begin{scope}[very thick]

\draw[->,brown] (2) -- (21);
\draw[->,brown] (2) -- (21');
\draw[->,blue] (2) -- (22);
\draw[dash pattern=on 2 off 1,brown] (21) -- (21');
\draw[dots,blue] (21) -- ++(-.7,.7);
\draw[dots,blue] (21) -- ++(.7,.7);
\draw[dots,blue] (21') -- ++(-.7,.7);
\draw[dots,blue] (21') -- ++(.7,.7);
\draw[dots,brown] (22) -- ++(-.9,.7) ++(.2,.05) coordinate (*);
\draw[dash pattern=on 2 off 1,brown,semithick] (*) -- ++(.6,0);
\draw[dots,brown] (22) -- ++(0,.7);
\draw[dots,blue] (22) -- ++(.9,.7);
\end{scope}
\end{scope}
\end{tikzpicture}
\caption{First levels of each cone type.}
\label{fig:different-cone-types-2}
\end{subfigure}
\hfil
\begin{subfigure}[t]{0.3\textwidth}
\centering
\begin{tikzpicture}[ very thick, every node/.style = {draw, circle, scale=.618}]
\node (0) at (0,0) {$0$};
\node (1) at (-1,1) {$1$};
\node (2) at (1,1) {$2$};

\begin{scope}[on background layer]
\begin{scope}[very thick]

\draw[->,brown] (0) to [out=175,in=-90,looseness=0.618] (1);
\draw[->,brown] (0) to [out=145,in=-60] (1);
\draw[->,blue] (0) to [out=35,in=240] (2);
\draw[->,blue] (0) to [out=5,in=270,looseness=0.618] (2);

\draw[->,blue] (1) to [out=-30,in=210,looseness=0.5] (2);
\draw[->,blue] (1) to [out=140,in=220,looseness=4] (1);
\draw[dash pattern=on 2 off 1,brown] (1) to [out=50,in=130,looseness=4] (1);

\draw[->,brown] (2) to [out=180,in=0,looseness=0.25] (1);
\draw[->,brown] (2) to [out=150,in=30,looseness=0.5] (1);
\draw[->,blue] (2) to [out=40,in=-40,looseness=4] (2);
\end{scope}
\end{scope}
\end{tikzpicture}
\caption{Graph of cone types.}
\label{fig:cone-types-2}
\end{subfigure}
\caption{Brown arrows represent right multiplication by $t$ or $\bar t$ while blue, by $u$ or $\bar u$. Dashed lines are length preserving.}
\end{figure}

\subsection*{Cayley graph}

\begin{figure}
\includegraphics[width=.618\textwidth]{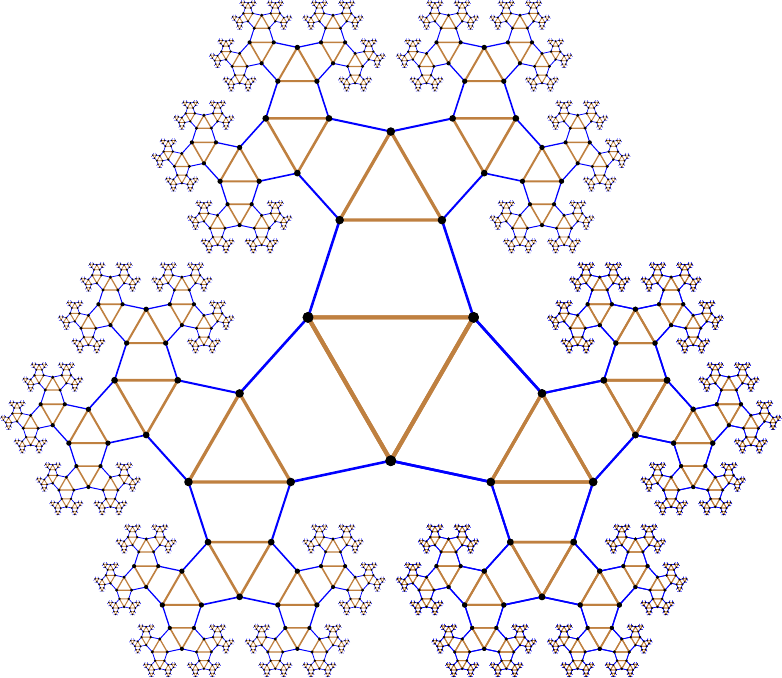}
\caption{The (undirected colored) Cayley graph $\Gamma(G,S)$ of $G = \pslz$ relative to $S'' = \{t,u\}$.
The $t$-edges in brown and the $u$-edges in blue.}
\label{fig:Cayley-graph-2}
\end{figure}

In the case of $G = \pslz$ and $S'' = \{t, u\}$, the Cayley graph $\Gamma(G,S'')$ is also weakly dual to a regular tree $\mathcal{T}$ of degree three, where the vertices in $\mathcal{T}$ correspond to triangles ($3$-cycles) in $\Gamma(G,S'')$ and the edges of $\mathcal{T}$ correspond to squares ($4$-cycles) in $\Gamma(G,S'')$ with opposite sides connecting two triangles. See \Cref{fig:Cayley-graph-2} for a representation of the (undirected colored) Cayley graph $\Gamma(G,S)$.
Each triangle corresponds to the primitive relators of length three: $t^2, \bar t^3$. The squares, to the primitive relators of length four: $(t\bar u)^2,(\bar tu)^2, (u\bar t)^2, (\bar ut)^2$.

In particular, each group element belongs to exactly one triangle and every triangle contains either one or two vertices of minimal length.
Moreover, this determines the cone types.
In fact, the cone type ---as $\typefunc''$ in \Cref{thm:cone-types-2}--- of an element $g \in G\setminus\{\id\}$ which is of minimal length in its triangle is $\typefunc''(g) = 2$ and otherwise $\typefunc''(g) = 1$.

Reciprocally, this geometric description of the cone types can be recovered from the combinatorial description: if $g \in G\setminus\{\id\}$, then $\typefunc''(g) = 1$ if and only if either $t$ or $\bar t$ is a suffix of some geodesic representing $g$, that is,  $\{t,\bar t\} \cap \suffix(g) \neq \emptyset$.

\subsection*{Asymptotic invariants}
We compute the growth, give upper bounds for the spectral radius and estimates for the entropy and the drift for $G = \pslz$ relative to $S'' = \{t,u\}$.

\begin{theorem}[Growth]
\label{thm:growth-2}
Let $G=\pslz$ and $S'' = \{t,u\}$.
The growth series of $G$ relative to $S''$ corresponds to the rational analytic function
\[
\frac{(1+z)^2}{1-2 z-z^2}.
\]
In particular, the rate of exponential growth of $G$ relative to $S''$ is
\[
v(G,S) = \sqrt{2} + 1
\approx
2.41 421
.\]
\end{theorem}

\begin{proof}
Let $g$ be the generating function for the spherical growth sequence $a_n = |\sphere_n|$, $n \in \N$.
If $g_t$ is the generating function for the spherical growth sequence of elements of $\typefunc''$-type $t \in \{0,1,2\}$, that is, of the sequence $a^{(t)}_n =  |\typefunc''^{-1}(t) \cap \sphere_n|$, $n \in \N$.

Now, \Cref{thm:cone-types-2} (cf. \Cref{fig:cone-types-2}) implies that there is the following recurrence relation:
\begin{gather*}
a^{(0)}_{0} = 1,
\qquad\text{and}\qquad
a^{(t)}_{0} = 0, \text{ for } t = 1,2,
\\[1ex]
a^{(0)}_{n+1} = 0,
\qquad
a^{(1)}_{n+1} = 2 a^{(0)}_{n} + a^{(1)}_n + 2 a^{(2)}_n,
\qquad
a^{(2)}_{n+1} = 2 a^{(0)}_{n} + a^{(1)}_{n} + a^{(2)}_n
\end{gather*}
and $a_{n} = a^{(0)}_{n} + a^{(1)}_{n} + a^{(2)}_{n}$.
Thus, the corresponding generating functions satisfy the equations
\[
g_0(z) = 1,
\qquad
g_1(z) = z (2 g_0(z) + g_1(z) + 2 g_2(z)),
\qquad
g_2(z) = z (2 g_0(z) + g_1(z) + g_2(z))
\]
and $g(z) = g_0(z) + g_1(z) + g_2(z)$.

Finally, solving the recurrence, we get that
\begin{gather*}
g_0(z) = 1,
\qquad
g_1(z) = \frac{2(1+z)z}{1-2 z - z^2},
\qquad
g_2(z) = \frac{2 z}{1-2 z - z^2}
\end{gather*}
and \[g(z) = \frac{(1+z)^2}{1-2 z-z^2}.\]

The growth rate $v(G,S)$ is the reciprocal of radius of convergence of $g$ around the origin, which corresponds to the smallest absolute value of a root of $1-2z-z^2$.
It follows that
\[
v(G,S) = \sqrt{2} + 1.
\qedhere
\]
\end{proof}


\begin{theorem}[Upper bound for the spectral radius]
\label{thm:upper-bound-2}
Let $G=\pslz$ and $S'' = \{t,u\}$.
Then, the spectral radius of the random walk on $G$ relative to $S''$ satisfies
\[
\rho(G,S'') <
0.926 762
.
\]
\end{theorem}
\begin{proof}
By \Cref{thm:cone-types}, the $f_t$'s of Nagnibeda's \Cref{thm:Nagnibeda} are given by:
\begin{align*}
f_0(c) & = 2 c_1 + 2 c_2, &
f_1(c) & = c_1 + c_2 + 1/c_1 + 1, &
f_2(c) & = 2 c_1 + c_2 + 1/c_2.
\end{align*}
It follows that $\rho(G,S'') \geq \max_{t} f_t(c)/|S'' \cup \bar S''|$, for every $c=(c_1,c_2)\in \R_+^{2}$.
Thus, the problem can be reduced to find the optimal such bound.
This can be solved numerically: we get that $\bar c\in\R_+^2$ with
\[
\bar c_1 \approx 0.703 111 830, \qquad
\bar c_2 \approx 0.581 687 196
\]
is a (global) minimun for $\max_{t} f_t(c)$, and
$\max_{t} f_t(\bar c) < 3.707 047 889$.

Finally, since $|S'' \cup \bar S''|=4$, by Nagnibeda's \Cref{thm:Nagnibeda}, it follows that
\[
\rho(G,S'') < 0.926 762.
\qedhere
\]
\end{proof}

\begin{remark}
In the case of $S'' = \{t, u\}$, since there is the relator $t^3$ of odd length, Gouëzel's \Cref{thm:Gouezel} does not apply.
In particular, we do not provide the analogue to \Cref{thm:lower-bound}.
Since the combinatorial description in this case is relatively simple, it is likely that estimates from below for $\rho(G,S'')$ can be computed by combinatorial methods.
We refrain from doing so here, however.
On the other hand, recall that we have the trivial lower bound given by the regular three of degree four.
Thus, we have
\[
\rho(G,S'') \geq \frac{\sqrt{3}}{2} \approx
0.866 025
.\]
\end{remark}

\begin{remark}
It is also worth to mention that, since in this case the Cayley graph is not bipartite (there are relators of odd length), the Markov spectrum is not necessarily symmetric. It follows that \cref{form:bottom-radius} is not valid, and we only have the inequality $\mu_0(G,S'') \leq 4(1-\rho(G,S''))$.
In particular, the estimates from above for the spectral radius $\rho(G,S'')$ as in \Cref{thm:upper-bound-2}, are of not use to estimate the bottom of the Laplace spectrum $\mu_0(G,S')$.
However, we still have the trivial upper bound given by the regular three of degree four.
Thus, we have
\[
\mu_0(G,S'') \geq 4-2\sqrt{3} \approx
0.535 898
.\]

\end{remark}

\begin{theorem}[Entropy and drift]
\label{thm:entropy-drift-2}
Let $G=\pslz$ and $S = \{r,u\}$.
Then, the entropy and the drift of the random walk on $G$ relative to $S$ satisfy
\[
0.296 750
< h(G,S) <
0.906 897
\qquad\text{and}\qquad
0.122 918
< \ell(G,S) <
0.375 650
.
\]
\end{theorem}

\begin{proof}
It is completely analogous to the proof of \Cref{thm:entropy-drift} in \Cref{sect:entropy-drift}, using the estimates from above of $\rho(G,S'')$ from \Cref{thm:upper-bound-2} and the value of $v(G,S'')$ from \Cref{thm:growth-2}.
\end{proof}

\section{Lower bounds for the spectral radius using Gouëzel's suffix types}
\label{sect:suffix-type}
In order to obtain better lower bounds for the spectral radius, we consider Gouëzel suffix-types introduced in \cite[Section~3.2]{Gouezel}.

Let $\tau\colon G \to T$ be any type function (or type system).
Then, for $g \in G$, consider the longest suffix that is common to every element in $\geodesics(g)$, say $w = w_1\dots w_n \in S^*$. Finally, we define the \defn{Gouëzel suffix type} of $g$ to be
\[\Gouezel(g) = \tau(g)\, \tau(g\bar w_n)\, \dots\, \tau(g\bar w_n \cdots \bar w_1) \in T^*.\]

In the case of the modular group $\pslz$ with generating system $S = \{r, u\}$, $\Gouezel$ is easy to compute inductively for $\tau = \typefunc$ (c.f \cite[Section~3.2]{Gouezel}).
In fact, we have that
\begin{itemize}
\item If $\typefunc(g) = 0$, that is, if $g = \id$, then $\Gouezel(\id)=0$.
\item If $\typefunc(g) = 5$, it has two predecessors, so $w$ is the empty word and $\Gouezel(g)=5$.
\item If $\typefunc(g) \in \{1, \dots, 4\}$, then $g$ has a unique predecessor, say $gs$; that is, $\{s\} = S_-(g)$.
Then, $\Gouezel(g)=\typefunc(g)\Gouezel(gs)$.
\end{itemize}

It also follows from this description that, if one knows $\Gouezel(g)$, one can determine inductively $\Gouezel(gs)$ for any $s \in S_+(g)$.
In fact, if $\typefunc(gs)=5$, then $\Gouezel(y)=5$, otherwise $g$ is the only predecessor of $gs$ and $\Gouezel(gs)=\typefunc(gs)\,\Gouezel(g)$.

In particular, $\Gouezel$ is a type function. However, these suffix types are not bounded in length and therefore do not define type systems as in \Cref{sect:Gouezel}. So it is natural to truncate the suffix types.
Gouëzel proposes several ways to do this in \cite[Section~3]{Gouezel}.
We follow the simplest, in which we fix a maximal length $n \in \N$, and define the \defn{Gouëzel suffix types of level~$n$} or \defn{$n$-suffix type} $\Gouezel_n(g)$ by keeping only the first $n$ elements of $\Gouezel(x)$.

It follows that $\Gouezel_n$ is a finite type function and, in particular, a type system.

It is clear that $\Gouezel_1$ coincides with $\typefunc$.
We study $\Gouezel_2$ in details in \Cref{sect:suffix-2}.
In \Cref{sect:suffix-higher} we summarize the corresponding results for $\Gouezel_n$ for $n = 3,\dots 12$.

\subsection{Estimates using Gouëzel suffix types of level \texorpdfstring{$2$}{2}}
\label{sect:suffix-2}
Following the previous discussion, we can compute all possible values for $\Gouezel_2$ and their corresponding successors.
See \Cref{tab:diagram-2}.

\begin{table}
\def\blanc{\phantom{0}}
\def\a{\texttt{a}}
\def\b{\texttt{b}}
\begin{tabular}{@{\quad} >{$}l<{$} @{$\quad\to\quad$} >{$}l<{$} @{\quad} >{$}r<{$} @{\,} >{$}l<{$} @{\quad}}
\toprule
\Gouezel_2(g) & \Gouezel_2(gs)\colon s\in S_+(g) \\
\midrule
0\blanc & 10, 20, 20 \\
10 & 31, 31 \\
2\a & 22, 42 & \text{for } \a & = 0,2,3 \\
3\b & 23, 5\blanc & \text{for } \b & = 1,4,5 \\
42 & 34, 5\blanc \\
5\blanc & 35 \\
\bottomrule
\end{tabular}
\caption{$\Gouezel_2$-types and the $\Gouezel_2$-types of their successors.}
\label{tab:diagram-2}
\end{table}

Types $0$, $10$, $20$ and $31$ represent a finite set of words of length at most two, so we can discard them to get a type system with $6 \times 6$ Perron-Frobenius matrices
\[
M = \begin{psmallmatrix}
 1 & 1 & 0 & 0 & 0 & 0 \\
 0 & 0 & 1 & 1 & 0 & 0 \\
 0 & 0 & 0 & 0 & 1 & 0 \\
 0 & 0 & 0 & 0 & 0 & 1 \\
 1 & 1 & 0 & 0 & 0 & 0 \\
 0 & 0 & 1 & 1 & 1 & 0
\end{psmallmatrix}
\qquad\text{and}\qquad
\tilde M = \begin{psmallmatrix}
 1 & 1 & 0 & 0 & 0 & 0 \\
 0 & 0 & 1 & 1 & 0 & 0 \\
 0 & 0 & 0 & 0 & 1 & 0 \\
 0 & 0 & 0 & 0 & 0 & 1 \\
 1 & 1 & 0 & 0 & 0 & 0 \\
 0 & 0 & \frac{1}{2} & \frac{1}{2} & \frac{1}{2} & 0
\end{psmallmatrix}.
\]
Then, we can proceed as in \Cref{sect:Gouezel} to obtain a lower bound for $\rho(G,S)$ using Gouëzel's \Cref{thm:Gouezel}.
For $\tilde M$, the Perron--Frobenius eigenvalue is $\eta = \frac{\sqrt{5}+1}{2}$ with an associated positive eigenvector
$
v = \left(
\frac{\sqrt{5}+1}{2}, \;
1, \;
1, \;
\frac{\sqrt{5}-1}{2}, \;
\frac{\sqrt{5}+1}{2}, \;
1
\right)^{T}
$.
Then, we get $D = \diag(v)$, $\hat M = D^{-1/2} M D^{1/2}$ and
\[
\bar M = \frac{\hat M + \hat M^{T}}{2} =
\frac{1}{2}
\begin{psmallmatrix}
 2 & \sqrt{\frac{\sqrt{5}-1}{2}} & 0 & 0 & 1 & 0 \\
 \sqrt{\frac{\sqrt{5}-1}{2}} & 0 & 1 & \sqrt{\frac{\sqrt{5}-1}{2}} & \sqrt{\frac{\sqrt{5}-1}{2}} & 0 \\
 0 & 1 & 0 & 0 & \sqrt{\frac{\sqrt{5}+1}{2}} & 1 \\
 0 & \sqrt{\frac{\sqrt{5}-1}{2}} & 0 & 0 & 0 & \sqrt{\frac{\sqrt{5}+1}{2}}^3 \\
 1 & \sqrt{\frac{\sqrt{5}-1}{2}} & \sqrt{\frac{\sqrt{5}+1}{2}} & 0 & 0 & \sqrt{\frac{\sqrt{5}+1}{2}} \\
 0 & 0 & 1 & \sqrt{\frac{\sqrt{5}+1}{2}}^3 & \sqrt{\frac{\sqrt{5}+1}{2}} & 0 
\end{psmallmatrix}.
\]

Finally, we compute numerically the dominant eigenvalue of the symmetric matrix $\bar M$, which is
$\lambda > 1.861 447 940 698$.
Thus, by Gouëzel's \Cref{thm:Gouezel}, we get that
\[
\rho(G,S) \geq \frac{2\lambda}{3\sqrt{\eta}} > 0.975 586.
\]

\subsection{Gouëzel suffix types of higher level}
\label{sect:suffix-higher}
Similar to \Cref{sect:suffix-2}, we can perform the analysis for $\Gouezel_n$, $n \geq 3$.
As final explicit examples, we include \Cref{tab:diagram-3,tab:diagram-4} showing all possible values for $\Gouezel_3$ and $\Gouezel_4$, respectively, with their corresponding successors.

\begin{table}
\def\blanc{\phantom{0}}
\def\bblanc{\phantom{00}}
\def\a{\texttt{a}}
\def\b{\texttt{b}}
\def\c{\texttt{c}}
\def\d{\texttt{d}}
\def\e{\texttt{e}}
\def\f{\texttt{f}}
\begin{tabular}{@{\quad} >{$}l<{$} @{$\quad\to\quad$} >{$}l<{$} @{\quad} >{$}r<{$} @{\,} >{$}l<{$} @{\quad}}
\toprule
\Gouezel_3(g) & \Gouezel_3(gs)\colon s\in S_+(g) \\
\midrule
0\bblanc & 10,\blanc 20,\blanc 20\blanc \\
10\blanc & 310, 310 \\
2\a\b & 22\a, 42\a & \text{for } \a\b & = 0\varepsilonup,20,22,23,31,34,35 \\
3\c\d & 23\c, 5\bblanc & \text{for } \c\d & = 10,42,5\varepsilonup \\
4\e\f & 34\e, 5\bblanc & \text{for } \e\f & = 20,22,23 \\
5\bblanc & 35\blanc \\
\bottomrule
\end{tabular}
\caption{$\Gouezel_3$-types and the $\Gouezel_3$-type of their successors.}
\label{tab:diagram-3}
\end{table}

\begin{table}
\def\blanc{\phantom{0}}
\def\bblanc{\phantom{00}}
\def\bbblanc{\phantom{000}}
\def\a{\texttt{a}}
\def\b{\texttt{b}}
\def\c{\texttt{c}}
\def\d{\texttt{d}}
\def\e{\texttt{e}}
\def\f{\texttt{f}}
\def\g{\texttt{g}}
\def\h{\texttt{h}}
\def\i{\texttt{i}}
\begin{tabular}{@{\quad} >{$}l<{$} @{$\quad\to\quad$} >{$}l<{$} @{\quad} >{$}r<{$} @{\,} >{$}l<{$} @{\quad}}
\toprule
\Gouezel_4(g) & \Gouezel_4(gs)\colon s\in S_+(g) \\
\midrule
0\bbblanc & 10,\bblanc 20,\bblanc 20\bblanc \\
10\bblanc & 310,\blanc 310\blanc \\
2\a\b\c & 22\a\b, 42\a\b & \text{for } \a\b\c & = 0\varepsilonup\varepsilonup,20\varepsilonup,222,223,231,234,235,310,342,35\varepsilonup \\
3\d\e\f & 23\d\e, 5\bbblanc & \text{for } \d\e\f & = 10\varepsilonup,420,422,423,5\varepsilonup\varepsilonup \\
4\g\h\i & 34\g\h, 5\bbblanc & \text{for } \g\h\i & = 20\varepsilonup,220,222,223,231,234,235,231,234,235 \\
5\bbblanc & 35\bblanc \\
\bottomrule
\end{tabular}
\caption{$\Gouezel_4$-types and the $\Gouezel_4$-type of their successors.}
\label{tab:diagram-4}
\end{table}

In the case of $\Gouezel_3$ one gets a type system with $\tilde{t}_3 = 9$ types (after excluding the finite set of elements of length at most $n=3$),
with $9 \times 9$ Perron-Frobenius matrices
\[
M = \begin{psmallmatrix}
 1 & 1 & 0 & 0 & 0 & 0 & 0 & 0 & 0 \\
 0 & 0 & 1 & 1 & 0 & 0 & 0 & 0 & 0 \\
 0 & 0 & 0 & 0 & 1 & 0 & 0 & 0 & 0 \\
 0 & 0 & 0 & 0 & 0 & 1 & 0 & 0 & 0 \\
 0 & 0 & 0 & 0 & 0 & 0 & 1 & 1 & 0 \\
 0 & 0 & 0 & 0 & 0 & 0 & 0 & 0 & 1 \\
 1 & 1 & 0 & 0 & 0 & 0 & 0 & 0 & 0 \\
 0 & 0 & 1 & 1 & 0 & 0 & 0 & 0 & 0 \\
 0 & 0 & 0 & 0 & 1 & 1 & 1 & 1 & 0 
\end{psmallmatrix}
\qquad\text{and}\qquad
\tilde M = \begin{psmallmatrix}
 1 & 1 & 0 & 0 & 0 & 0 & 0 & 0 & 0 \\
 0 & 0 & 1 & 1 & 0 & 0 & 0 & 0 & 0 \\
 0 & 0 & 0 & 0 & 1 & 0 & 0 & 0 & 0 \\
 0 & 0 & 0 & 0 & 0 & 1 & 0 & 0 & 0 \\
 0 & 0 & 0 & 0 & 0 & 0 & 1 & 1 & 0 \\
 0 & 0 & 0 & 0 & 0 & 0 & 0 & 0 & 1 \\
 1 & 1 & 0 & 0 & 0 & 0 & 0 & 0 & 0 \\
 0 & 0 & 1 & 1 & 0 & 0 & 0 & 0 & 0 \\
 0 & 0 & 0 & 0 &  \frac{1}{2} & \frac{1}{2} & \frac{1}{2} & \frac{1}{2} & 0 
\end{psmallmatrix},
\]
and the estimate
$
\rho(G,S) \geq \tilde\rho_3 =
0.975 680
$.
For $\Gouezel_4$ one gets $\tilde{t}_4 = 14$ types and 
$
\tilde\rho_4 = 0.975 712
$.
For higher levels, up to $n = 12$, see \Cref{tab:lower-bounds}.

\begin{table}
\begin{tabular}{@{\quad} c @{\qquad} c @{\qquad} c @{\quad}}
\toprule
Level & Number of types & Lower bound \\
\midrule
$2$ & $6$ & $0.975\, 586\, 575$ \\
$3$ & $9$ & $0.975\, 680\, 157$ \\
$4$ & $14$ & $0.975\, 712\, 971$ \\
$5$ & $21$ & $0.975\, 723\, 071$ \\
$6$ & $31$ & $0.975\, 726\, 073$ \\
$7$ & $46$ & $0.975\, 727\, 002$ \\
$8$ & $68$ & $0.975\, 727\, 290$ \\
$9$ & $100$ & $0.975\, 727\, 378$ \\
$10$ & $147$ & $0.975\, 727\, 405$ \\
$11$ & $216$ & $0.975\, 727\, 414$ \\
$12$ & $317$ & $0.975\, 727\, 416$ \\
\bottomrule
\end{tabular}%
\medskip
\caption{Lower bounds for the spectral radius using Gouëzel $n$-suffix types.}
\label{tab:lower-bounds}
\end{table}

\begin{remark}
The main issue in using Gouëzel suffix types is that they do not distinguish group elements with several predecessors.
In our case, these correspond to type-$5$ elements.
In \cite[Section~3.3]{Gouezel}, Gouëzel introduce \emph{essential types} that allows to partially overcome this issue.
In \Cref{sect:my-type}, we consider a new kind of suffix types that allows better numerical estimates.
These new suffix types share the essence of Gouëzel suffix and essential types in \cite[Section~3]{Gouezel} and also of Cannon's $N$-types in \cite[Section~7]{Cannon}.
\end{remark}

\section{Lower bounds using \defn{new} suffix types}
\label{sect:my-type}
In this appendix, we introduce another kind of suffix types that allows better numerical estimates from below for the spectral radius than those obtained in \Cref{sect:suffix-type}.

By \Cref{prop:sufre-type}, we know that if the number of primitive relators is finite, then, it is enough to know long enough suffixes of geodesics in order to determine the cone type of an element.
This motivates the following suffix type functions.

Recall that $\suffix(g)$ is the set of all suffixes of geodesics for $g\in G$.
The \defn{$n$-suffix} of $g$, denoted $\suffix_n(g)$, is the set of all maximal suffixes of length at most $n\in \N$ of geodesics for $g\in G$, that is,
\[\suffix_n(g) = \suffix(g) \cap \sphere_{\min\{n,|g|_S\}}.\]

In the case when the set of primitive relators is finite --as is the case of $\pslz$--, the type function $\sufre$ is finite and therefore, for $n$ large enough, $\suffix_n$ is also a type function.
In particular, in our case of interest, we have the following.

\begin{theorem}
\label{thm:my-type}
For any $n \geq 2$, the $n$-suffix of an element in $G = \pslz$ relative to $S = \{r,u\}$ determines its cone type.
That is, $\suffix_n\colon G \to 2^{S^*}$ is a type function.
\end{theorem}
\begin{proof}
Since $\suffix_{n+1}$ determines $\suffix_n$, it is enough to prove the result for $n=2$.
Before doing so, let us notice first that, since every element in $\sufre(g)$ is a suffix of $g$ of length at most three, it is clear that $\suffix_3(g)$ determines $\sufre(g)$ and therefore, by \Cref{prop:sufre-type}, its cone type.

Now, in order to conclude in the case of $\suffix_2$, it is enough to note that the only $\sufre$-types that have suffixes of length greater than two are $\{rur,\bar ur\bar u\}$ and $\{r\bar ur,uru\}$. But, in these cases, the only possibility is for $\suffix_2$ to be $\{r\bar u,ur\}$ and $\{ru,\bar ur\}$, respectively.
Thus, $\suffix_2$ determines $\sufre(g)$. Again, by \Cref{prop:sufre-type}, we conclude that $\suffix_2$ determines the cone types.
\end{proof}

\begin{table}
\begin{tabular}{@{\quad} >{$}l<{$} @{$\quad\to\quad$} >{$}l<{$} @{\quad}}
\toprule
\suffix_2(g) & \suffix_2(gs)\colon s\in S_+(g) \\
\midrule
\emptyset & \{r\}, \{u\}, \{\bar u\} \\
\{r\} & \{ru\}, \{r\bar u\} \\
\{a\} \text{ or } \{a^2\} & \{ar\}, \{a^2\} \\
\{ra\} \text{ or } \{ra, a^2\} & \{a^2\}, \{r\bar a,ar\}\\
\{ar\} & \{r\bar a\}, \{ra,\bar ar\} \\
\{ra,\bar ar\} & \{ra, a^2\} \\
\bottomrule
\end{tabular}
\caption{$\suffix_2$-types and the $\suffix_2$-type of their successors. Here, $a = u$ or $\bar u$.}
\label{tab:diagram-my2}
\end{table}

\subsection{Computing higher level suffix types}
\label{sect:recursion}
As discussed above, when the number of primitive relators is finite, for $n$ large enough, $\suffix_n$ completely determines $\sufre$ and, in particular, the cone types.
Moreover, it is also possible to determine $\suffix_{n+1}$ inductively.

In the particular case of $G = \pslz$ and $S = \{r,u\}$, one has that $s \in S_+(g)$ if and only if $s \notin \suffix(g)$.
In particular, we can determine $S_+(g)$ from $\suffix_n(g)$ for any $n > 0$.
Moreover, if $n \geq 2$, from the description of the primitive relators, one has the following rules to determine $\suffix_{n+1}$ from $\suffix_n$ (cf. \Cref{lemm:successor-type}):
\begin{itemize}
\item If $s \in S_+(g)$, then $ws \in \suffix_{n+1}(gs)$ for each $w \in \suffix_n(g)$;
\item For $a \in \{u,\bar u\}$, if $a \in S_+(g)$, then $wr\bar ar \in \suffix_{n+1}(ga)$ for each $war \in \suffix_n(g)$; and
\item If $r \in S_+(g)$, then $wara \in \suffix_{n+1}(gr)$ for each $wr\bar a \in \suffix_n(g)$, $a \in \{u,\bar u\}$.
\end{itemize}
Similarly, one can compute the $n$-suffix of the successors of a given $n$-suffix, when $n \geq 2$.

\subsection{Estimates from suffix types of level \texorpdfstring{$2$}{2}}
\label{sect:type-2}
In \Cref{tab:diagram-my2}, we summarize the $\suffix_2$-types together with their corresponding successors.
There is no gain in differentiating $\suffix_2$-types that coincide after exchanging $u$ with $\bar u$. So we consider them as the same type.
Note also that if we identify the $\suffix_2$-types that have the same successors' $\suffix_2$-types, then $\suffix_2$ is equivalent to $\sufre$ as type function, giving an alternative proof for \Cref{thm:my-type}.
This also tell us that we cannot get better estimates with these latter identifications.
If we do not identify them, then we get slightly better estimates.

In fact, by discarding the finite set of words of length one, we get a type system with $5 \times 5$ Perron-Frobenius matrices $M$ and $\tilde M$, where
\[
M = \begin{psmallmatrix}
 1 & 1 & 1 & 0 & 0 \\
 0 & 0 & 0 & 1 & 0 \\
 0 & 0 & 0 & 0 & 1 \\
 1 & 0 & 0 & 0 & 0 \\
 0 & 1 & 1 & 1 & 0
\end{psmallmatrix}
\qquad\text{and}\qquad
\tilde M = \begin{psmallmatrix}
 1 & 1 & 1 & 0 & 0 \\
 0 & 0 & 0 & 1 & 0 \\
 0 & 0 & 0 & 0 & 1 \\
 1 & 0 & 0 & 0 & 0 \\
 0 & \frac{1}{2} & \frac{1}{2} & \frac{1}{2} & 0
\end{psmallmatrix}.
\]
Following Gouëzel's \Cref{thm:Gouezel}, for $\tilde M$, we compute the Perron--Frobenius eigenvalue $\eta = \frac{\sqrt{5}+1}{2}$ and an associated positive eigenvector
$
v = \left(
\frac{\sqrt{5}+3}{2}, \;
1, \;
\frac{\sqrt{5}-1}{2}, \;
\frac{\sqrt{5}+1}{2}, \;
1
\right)^{T}
$.
Then, we get $D = \diag(v)$, $\hat M = D^{-1/2} M D^{1/2}$ and
\[
\bar M = \frac{\hat M + \hat M^{T}}{2} =
\frac{1}{2}
\begin{psmallmatrix}
 2 & \frac{\sqrt{5}-1}{2} & \sqrt{\frac{\sqrt{5}-1}{2}}^3 & \sqrt{\frac{\sqrt{5}+1}{2}} & 0 \\
 \frac{\sqrt{5}-1}{2} & 0 & 0 & \sqrt{\frac{\sqrt{5}+1}{2}} & 1 \\
 \sqrt{\frac{\sqrt{5}-1}{2}}^3 & 0 & 0 & 0 & \sqrt{\frac{\sqrt{5}+1}{2}}^3 \\
 \sqrt{\frac{\sqrt{5}+1}{2}} & \sqrt{\frac{\sqrt{5}+1}{2}} & 0 & 0 & \sqrt{\frac{\sqrt{5}+1}{2}} \\
 0 & 1 & \sqrt{\frac{\sqrt{5}+1}{2}}^3 & \sqrt{\frac{\sqrt{5}+1}{2}} & 0
\end{psmallmatrix}.
\]

Finally, the dominant eigenvalue of $\bar M$ satisfies
$\lambda > 1.861 199 191 510$
and, by Gouëzel's \Cref{thm:Gouezel},
\[
\rho(G,S) \geq \frac{2\lambda}{3\sqrt{\eta}} > 0.975 456,
\]
slightly improving \Cref{thm:lower-bound}, but still worse than using Gouëzel suffix types of level~$2$ as in \Cref{sect:suffix-2}.

\subsection{Suffix types of higher level}
\label{sect:type-higher}
Similar to \Cref{sect:type-2}, and following \Cref{sect:recursion}, we can perform the analysis for $\suffix_n$, $n \geq 3$.
As explicit examples, in \Cref{tab:diagram-my3,tab:diagram-my4} we include the type systems corresponding to $\suffix_3$ and $\suffix_4$, respectively, with their corresponding successors.

\begin{table}
\begin{tabular}{@{\quad} >{$}l<{$} @{$\quad\to\quad$} >{$}l<{$} @{\quad}}
\toprule
\suffix_3(g) & \suffix_3(gs)\colon s\in S_+(g) \\
\midrule
\{ra^2\}, \; \{a^3\} \text{ or } \{ra^2,a^3\}  & \{a^2r\}, \{a^3\} \\
\{ar\bar a\} & \{r\bar a^2\}, \{r\bar ar,ara\} \\
\{a^2r\} & \{ar\bar a\}, \{r\bar ar,ara\} \\
\{rar,\bar a r\bar a\} \text{ or } \{rar,\bar a r\bar a,a^2r\} & \{r\bar a^2,ar\bar a\} \\
\{ra^2,\bar ara\} & \{ra^2,a^3\}, \{ra r,\bar ar\bar a,a^2r\} \\
\bottomrule
\end{tabular}
\caption{$\suffix_3$-types and the $\suffix_3$-type of their successors. Here, $a = u$ or $\bar u$.
This table excludes the $\suffix_3$-types of elements of length bounded by three.}
\label{tab:diagram-my3}
\end{table}

\begin{table}
\begin{tabular}{@{\quad} >{$}l<{$} @{$\quad\to\quad$} >{$}l<{$} @{\quad}}
\toprule
\suffix_4(g) & \suffix_4(gs)\colon s\in S_+(g) \\
\midrule
\{ra^2r\}, \; \{a^3r\}, \; \{ra^2r,a^3r\} & \{a^2r\bar a\}, \{ar\bar ar,a^2ra\} \\
\{ra^3\}, \; \{a^4\}, \; \{ra^3, a^4\} & \{a^3r\}, \{a^4\} \\
\{ar\bar a^2\} & \{r\bar a^2r\}, \{r\bar a^3\} \\
\{a^2r\bar a\} & \{ar\bar a^2\}, \{ar\bar a r,a^2ra\} \\
\{rar\bar a,\bar ar\bar a^2\}, \; \{rar\bar a,\bar ar\bar a^2,a^2r\bar a\} & \{r\bar a^3,ar\bar a^2\}, \{r\bar a^2r,ar\bar ar,a^2ra\} \\
\{ra^3,\bar ara^2\} & \{ra^2r, a^3r\}, \{ra^3,a^4\} \\
\{ar\bar ar,a^2ra\}  & \{r\bar ara,ara^2\} \\
\{ra^2r,\bar ara r,\bar a^2r\bar a\} & \{rar\bar a,\bar ar\bar a^2,a^2r\bar a\} \\
\bottomrule
\end{tabular}
\caption{$\suffix_4$-types and the $\suffix_4$-type of their successors.Here, $a = u$ or $\bar u$.
This table excludes the $\suffix_4$-types of elements of length bounded by four.}
\label{tab:diagram-my4}
\end{table}

In the case of $\suffix_3$ one gets a type system with $t_3 = 8$ types (after excluding the finite set of elements of length at most $n=3$)
with $8 \times 8$ Perron-Frobenius matrices
\[
M = \begin{psmallmatrix}
 0 & 1 & 0 & 0 & 0 & 0 & 0 & 0 \\
 0 & 0 & 1 & 1 & 1 & 0 & 0 & 0 \\
 1 & 0 & 0 & 0 & 0 & 0 & 0 & 0 \\
 0 & 0 & 1 & 1 & 1 & 0 & 0 & 0 \\
 0 & 0 & 0 & 0 & 0 & 0 & 0 & 1 \\
 1 & 1 & 0 & 0 & 0 & 0 & 0 & 0 \\
 0 & 0 & 0 & 0 & 0 & 0 & 0 & 1 \\
 0 & 0 & 0 & 0 & 0 & 1 & 1 & 0
\end{psmallmatrix}
\qquad\text{and}\qquad
\tilde M = \begin{psmallmatrix}
 0 & 1 & 0 & 0 & 0 & 0 & 0 & 0 \\
 0 & 0 & 1 & 1 & 1 & 0 & 0 & 0 \\
 1 & 0 & 0 & 0 & 0 & 0 & 0 & 0 \\
 0 & 0 & 1 & 1 & 1 & 0 & 0 & 0 \\
 0 & 0 & 0 & 0 & 0 & 0 & 0 & 1 \\
 \frac{1}{2} & \frac{1}{2} & 0 & 0 & 0 & 0 & 0 & 0 \\
 0 & 0 & 0 & 0 & 0 & 0 & 0 & \frac{1}{2} \\
 0 & 0 & 0 & 0 & 0 & 1 & 1 & 0
\end{psmallmatrix},
\]
and the estimate
$
\rho(G,S) \geq \rho_3 = 0.975 847
$.
For $\suffix_4$ one gets $t_4 = 13$ types and 
$
\rho_4 = 0.976 000
$.
For higher levels, up to $n = 10$, see \Cref{tab:my-lower-bounds}.
\smallbreak
This also concludes the proof of \Cref{thm:spectral-radius}.
\pushQED{\qed}\qedhere

\begin{table}
\begin{tabular}{@{\quad} c @{\qquad} c @{\qquad} c @{\quad}}
\toprule
Level & Number of types & Lower bound \\
\midrule
$2$ & $5$ & $0.975\, 456\, 205$ \\
$3$ & $8$ & $0.975\, 847\, 429$ \\
$4$ & $13$ & $0.976\, 000\, 849$ \\
$5$ & $23$ & $0.976\, 112\, 118$ \\
$6$ & $38$ & $0.976\, 183\, 437$ \\
$7$ & $64$ & $0.976\, 238\, 443$ \\
$8$ & $107$ & $0.976\, 278\, 755$ \\
$9$ & $181$ & $0.976\, 310\, 925$ \\
$10$ & $303$ & $0.976\, 336\, 575$ \\
\bottomrule
\end{tabular}%
\medskip
\caption{Lower bounds for the spectral radius using our $n$-suffix types.}
\label{tab:my-lower-bounds}
\end{table}

\printbibliography
\end{document}